\documentclass {article}
\usepackage{amsmath,amssymb}
\usepackage{cite}
\begin{document}

\title{Determinantal representations of the Drazin inverse for Hermitian matrix over the
quaternion skew field with applications.}
\author{Ivan Kyrchei \footnote{Pidstrygach Institute for Applied Problems of Mechanics and Mathematics of NAS of Ukraine,
str. Naukova 3b, Lviv, 79060, Ukraine, kyrchei@online.ua}}

\date{}
 \maketitle

\begin{abstract}
Within the framework of the theory of the column and row
determinants, we obtain determinantal representations of the Drazin inverse for Hermitian matrix over the
quaternion skew field. Using the obtained determinantal representations of the Drazin inverse we get explicit representation formulas (analogs
of Cramer's rule) for the Drazin inverse solutions of
quaternion matrix equations $ {\bf A}{\bf X} = {\bf D}$, $ {\bf X}{\bf B} = {\bf D} $ and ${\bf A} {\bf X} {\bf B} = {\bf D} $, where $
{\bf A}$, ${\bf B}$ are Hermitian.
\end{abstract}

\textit{Keyword}: Matrix equation,  Drazin inverse solution,
Drazin inverse, Quaternion matrix,
Cramer rule, Column determinant,  Row determinant.

{\bf MSC}: 15A15, 16W10.

\section{Introduction}
\newtheorem{Corollary}{Corollary}[section]
\newtheorem{theorem}{Theorem}[section]
\newtheorem{lemma}{Lemma}[section]
\newtheorem{definition}{Definition}[section]
\newtheorem{remark}{Remark}[section]
\newcommand{\rank}{\mathop{\rm rank}\nolimits}
\newtheorem{proposition}{Proposition}[section]

Throughout the paper, we denote the real number field by ${\rm
{\mathbb{R}}}$, the set of all $m\times n$ matrices over the
quaternion algebra
\[{\rm {\mathbb{H}}}=\{a_{0}+a_{1}i+a_{2}j+a_{3}k\,
|\,i^{2}=j^{2}=k^{2}=-1,\, a_{0}, a_{1}, a_{2}, a_{3}\in{\rm
{\mathbb{R}}}\}\]
by ${\rm {\mathbb{H}}}^{m\times n}$.  Let ${\rm M}\left( {n,{\rm {\mathbb{H}}}} \right)$ be the
ring of $n\times n$ quaternion matrices. For ${\rm {\bf A}}
 \in {\rm {\mathbb{H}}}^{n\times m}$, the symbols ${\rm {\bf A}}^{ *}$ stands for the conjugate transpose (Hermitian adjoint) matrix
of ${\rm {\bf A}}$.
 The matrix ${\rm {\bf A}} = \left( {a_{ij}}  \right) \in {\rm
{\mathbb{H}}}^{n\times n}$ is Hermitian if ${\rm {\bf
A}}^{ *}  = {\rm {\bf A}}$.

As one of the important types of generalized inverses of matrices, the Drazin inverses and their applications
have well been examined in the literature (see, e.g., \cite{dr,ca1,be,zh1,ha,kh}). Stanimirovi\`{c} and Djordjevi\'{c} in \cite{st} have introduced a determinantal representation of the Drazin
inverse of a complex matrix based on its full-rank representation. In \cite{ky1} we obtain  determinantal representations of the complex Moore-Penrose, and Drazin inverses, respectively,  used its limit representation. It allowed to obtain \cite{ky11,ky2}  the analogs of Cramer's rule for the least squares solutions with the minimum norm and the Drazin inverse solutions of the following matrix equations
\begin{equation}\label{eq12:AX=D}
{\rm {\bf A}}{\rm {\bf X}} = {\rm {\bf D}},
\end{equation}
\begin{equation}\label{eq12:XB=D}
{\rm {\bf X}}{\rm {\bf B}} ={\rm {\bf D}}, \end{equation}
\begin{equation}\label{eq12:AXB=D} {\rm
{\bf A}}{\rm {\bf X}}{\rm {\bf B}} ={\rm {\bf D}}.\end{equation}
In the case of quaternion matrices we are faced with the problem of the determinant of a square quaternion matrix.  But recently we have developed  the theory
of the column and row determinants of a quaternion matrix in \cite{ky3,ky4}.
Within the framework of the theory of the column and row determinants, we have obtained determinantal representations of  the Moore-Penrose inverse  by analogs of the classical adjoint matrix  in \cite{ky5,ky33}  and analogs of Cramer's rule for the least squares solutions with minimum norm of the  matrix equations (\ref{eq12:AX=D}), (\ref{eq12:XB=D}), and (\ref{eq12:AXB=D}) in \cite{ky6}.

In \cite{so2,so1,song}
the authors have received determinantal representations of the generalized
inverses ${\rm {\bf A}}^{2}_{r_{T_{1},S_{1}}}$,  ${\rm {\bf A}}^{2}_{l_{T_{2},S_{2}}}$,  and  ${\rm {\bf A}}^{2}_{(T_{1},T_{2}), (S_{1},S_{2})}$, and consequently of the Moore-Penrose and Drazin inverses
over the quaternion skew field by the theory of the
column and row determinants  as well. But in obtaining of these determinantal representations another auxiliary matrices together with ${\rm {\bf A}}$  are used.

In this paper we aim to obtain determinantal representations of the Drazin inverse for a Hermitian quaternion  matrix by using only entries of ${\rm {\bf A}}$ and
 explicit representation formulas (analogs
of Cramer's rule) for the Drazin inverse solutions of quaternion matrix equations  (\ref{eq12:AX=D}), (\ref{eq12:XB=D}), and (\ref{eq12:AXB=D}), where ${\rm {\bf A}}$ and ${\rm {\bf B}}$ are Hermitian, without any restriction. Obtaining of  determinantal representation of the Drazin inverse for an arbitrary quaternion matrix is a difficult task that requires more research.

 The paper is organized as follows. We start with some
basic concepts and results from the theory of  the row and column determinants and the theory on eigenvalues   of quaternion  matrices  in Section 2. We  give the determinantal representations of the Drazin inverse for a Hermitian  quaternion matrix  in Section 3. In Section 4, we obtain
explicit representation formulas for the Drazin inverse solutions of quaternion matrix equations  (\ref{eq12:AX=D}), (\ref{eq12:XB=D}), and (\ref{eq12:AXB=D}). In Section 5, we show
a numerical example to illustrate the main result.

\section{Elements of the theory of the column and row determinants.}
Suppose $S_{n}$ is the
symmetric group on the set $I_{n}=\{1,\ldots,n\}$.
\begin{definition}
 The $i$th row determinant of ${\rm {\bf A}}=(a_{ij}) \in {\rm
M}\left( {n,{\mathbb{H}}} \right)$ is defined by
 \[{\rm{rdet}}_{ i} {\rm {\bf A}} =
{\sum\limits_{\sigma \in S_{n}} {\left( { - 1} \right)^{n -
r}{a_{i{\kern 1pt} i_{k_{1}}} } {a_{i_{k_{1}}   i_{k_{1} + 1}}}
\ldots } } {a_{i_{k_{1} + l_{1}}
 i}}  \ldots  {a_{i_{k_{r}}  i_{k_{r} + 1}}}
\ldots  {a_{i_{k_{r} + l_{r}}  i_{k_{r}} }}\] \noindent for all $i
= 1,\ldots,n $. The left-ordered cycle notation of the
permutation $\sigma$ is written as follows,
\[\sigma = \left(
{i\,i_{k_{1}}  i_{k_{1} + 1} \ldots i_{k_{1} + l_{1}} }
\right)\left( {i_{k_{2}}  i_{k_{2} + 1} \ldots i_{k_{2} + l_{2}} }
\right)\ldots \left( {i_{k_{r}}  i_{k_{r} + 1} \ldots i_{k_{r} +
l_{r}} } \right).\] \noindent The index $i$ opens the first cycle
from the left  and other cycles satisfy the following conditions,
$i_{k_{2}} < i_{k_{3}}  < \ldots < i_{k_{r}}$ and $i_{k_{t}}  <
i_{k_{t} + s} $ for all $t = 2,\ldots,r $ and $s =1,\ldots,l_{t}
$.
\end{definition}

\begin{definition}
The $j$th column determinant
 of ${\rm {\bf
A}}=(a_{ij}) \in {\rm M}\left( {n,{\mathbb{H}}} \right)$ is
defined by
 \[{\rm{cdet}} _{{j}}\, {\rm {\bf A}} =
{{\sum\limits_{\tau \in S_{n}} {\left( { - 1} \right)^{n -
r}a_{j_{k_{r}} j_{k_{r} + l_{r}} } \ldots a_{j_{k_{r} + 1}
i_{k_{r}} }  \ldots } }a_{j\, j_{k_{1} + l_{1}} }  \ldots  a_{
j_{k_{1} + 1} j_{k_{1}} }a_{j_{k_{1}} j}}\] \noindent for all $j
=1,\ldots,n $. The right-ordered cycle notation of the permutation
$\tau \in S_{n}$ is written as follows,
 \[\tau =
\left( {j_{k_{r} + l_{r}}  \ldots j_{k_{r} + 1} j_{k_{r}} }
\right)\ldots \left( {j_{k_{2} + l_{2}}  \ldots j_{k_{2} + 1}
j_{k_{2}} } \right){\kern 1pt} \left( {j_{k_{1} + l_{1}}  \ldots
j_{k_{1} + 1} j_{k_{1} } j} \right).\] \noindent The index $j$
opens  the first cycle from the right  and other cycles satisfy
the following conditions, $j_{k_{2}}  < j_{k_{3}}  < \ldots <
j_{k_{r}} $ and $j_{k_{t}}  < j_{k_{t} + s} $ for all $t =
2,\ldots,r $ and $s = 1,\ldots,l_{t}  $.
\end{definition}

Suppose ${\rm {\bf A}}_{}^{i{\kern 1pt} j} $ denotes the submatrix
of ${\rm {\bf A}}$ obtained by deleting both the $i$th row and the
$j$th column. Let ${\rm {\bf a}}_{.j} $ be the $j$th column and
${\rm {\bf a}}_{i.} $ be the $i$th row of ${\rm {\bf A}}$. Suppose
${\rm {\bf A}}_{.j} \left( {{\rm {\bf b}}} \right)$ denotes the
matrix obtained from ${\rm {\bf A}}$ by replacing its $j$th column
with the column ${\rm {\bf b}}$, and ${\rm {\bf A}}_{i.} \left(
{{\rm {\bf b}}} \right)$ denotes the matrix obtained from ${\rm
{\bf A}}$ by replacing its $i$th row with the row ${\rm {\bf b}}$.

We  note some properties of column and row determinants of a
quaternion matrix ${\rm {\bf A}} = \left( {a_{ij}} \right)$, where
$i \in I_{n} $, $j \in J_{n} $ and $I_{n} = J_{n} = {\left\{
{1,\ldots ,n} \right\}}$.
\begin{proposition}  \cite{ky3}
If $b \in {\mathbb{H}}$, then
 $ {\rm{rdet}}_{ i} {\rm {\bf A}}_{i.} \left( {b
\cdot {\rm {\bf a}}_{i.}}  \right) = b \cdot {\rm{rdet}}_{ i} {\rm
{\bf A}}$ for all $i =1,\ldots,n $.
\end{proposition}
\begin{proposition}\label{pr:col_mult} \cite{ky3}
If $b \in {\mathbb{H}}$, then  ${\rm{cdet}} _{{j}}\, {\rm {\bf
A}}_{.j} \left( {{\rm {\bf a}}_{.j} \cdot b} \right) = {\rm{cdet}}
_{{j}}\, {\rm {\bf A}} \cdot  b$ for all $j =1,\ldots,n$.
\end{proposition}
\begin{proposition} \cite{ky3}
If for  ${\rm {\bf A}}\in {\rm M}\left( {n,{\mathbb{H}}}
\right)$\, there exists $t \in I_{n} $ such that $a_{tj} = b_{j} +
c_{j} $\, for all $j = 1,\ldots,n$, then
\[
\begin{array}{l}
   {\rm{rdet}}_{{i}}\, {\rm {\bf A}} = {\rm{rdet}}_{{i}}\, {\rm {\bf
A}}_{{t{\kern 1pt}.}} \left( {{\rm {\bf b}}} \right) +
{\rm{rdet}}_{{i}}\, {\rm {\bf A}}_{{t{\kern 1pt}.}} \left( {{\rm
{\bf c}}} \right), \\
  {\rm{cdet}} _{{i}}\, {\rm {\bf A}} = {\rm{cdet}} _{{i}}\, {\rm
{\bf A}}_{{t{\kern 1pt}.}} \left( {{\rm {\bf b}}} \right) +
{\rm{cdet}}_{{i}}\, {\rm {\bf A}}_{{t{\kern 1pt}.}} \left( {{\rm
{\bf c}}} \right),
\end{array}
\]
\noindent where ${\rm {\bf b}}=(b_{1},\ldots, b_{n})$, ${\rm {\bf
c}}=(c_{1},\ldots, c_{n})$ and for all $i =1,\ldots,n$.
\end{proposition}
\begin{proposition} \cite{ky3}
If for ${\rm {\bf A}}\in {\rm M}\left( {n,{\mathbb{H}}} \right)$\,
 there exists $t \in J_{n} $ such that $a_{i\,t} = b_{i} + c_{i}$
for all $i = 1,\ldots,n$, then
\[
\begin{array}{l}
  {\rm{rdet}}_{{j}}\, {\rm {\bf A}} = {\rm{rdet}}_{{j}}\, {\rm {\bf
A}}_{{\,.\,{\kern 1pt}t}} \left( {{\rm {\bf b}}} \right) +
{\rm{rdet}}_{{j}}\, {\rm {\bf A}}_{{\,.\,{\kern 1pt} t}} \left(
{{\rm
{\bf c}}} \right),\\
  {\rm{cdet}} _{{j}}\, {\rm {\bf A}} = {\rm{cdet}} _{{j}}\, {\rm
{\bf A}}_{{\,.\,{\kern 1pt}t}} \left( {{\rm {\bf b}}} \right) +
{\rm{cdet}} _{{j}} {\rm {\bf A}}_{{\,.\,{\kern 1pt}t}} \left(
{{\rm {\bf c}}} \right),
\end{array}
\]
\noindent where ${\rm {\bf b}}=(b_{1},\ldots, b_{n})^T$, ${\rm
{\bf c}}=(c_{1},\ldots, c_{n})^T$ and for all $j =1,\ldots,n$.
\end{proposition}
\begin{proposition}\cite{ky3}
If ${\rm {\bf A}}^{ *} $
is the Hermitian adjoint matrix of  ${\rm {\bf A}}\in {\rm
M}\left( {n,{\rm {\mathbb{H}}}} \right)$, then $
{\rm{rdet}}_{{i}}\, {\rm {\bf A}}^{ *} =  \overline{{{\rm{cdet}}
_{{i}}\, {\rm {\bf A}}}}$ for all $i = 1,\ldots,n $.
\end{proposition}
The following lemmas enable us to expand ${\rm{rdet}}_{{i}}\, {\rm
{\bf A}}$ by cofactors
  along  the $i$th row and ${\rm{cdet}} _{j} {\rm {\bf A}}$
 along  the $j$th column respectively for all $i, j = 1,\ldots,n$.

\begin{lemma}\label{lemma:R_ij} \cite{ky3}
Let $R_{i{\kern 1pt} j}$ be the right $ij$-th cofactor of ${\rm
{\bf A}}\in {\rm M}\left( {n, {\mathbb{H}}} \right)$, that is, $
{\rm{rdet}}_{{i}}\, {\rm {\bf A}} = {\sum\limits_{j = 1}^{n}
{{a_{i{\kern 1pt} j} \cdot R_{i{\kern 1pt} j} } }} $ for all $i =
1,\ldots,n$.  Then
\[
 R_{i{\kern 1pt} j} = {\left\{ {{\begin{array}{*{20}c}
  - {\rm{rdet}}_{{j}}\, {\rm {\bf A}}_{{.{\kern 1pt} j}}^{{i{\kern 1pt} i}} \left( {{\rm
{\bf a}}_{{.{\kern 1pt} {\kern 1pt} i}}}  \right),& {i \ne j},
\hfill \\
 {\rm{rdet}} _{{k}}\, {\rm {\bf A}}^{{i{\kern 1pt} i}},&{i = j},
\hfill \\
\end{array}} } \right.}
\]
\noindent where  ${\rm {\bf A}}_{.{\kern 1pt} j}^{i{\kern 1pt} i}
\left( {{\rm {\bf a}}_{.{\kern 1pt} {\kern 1pt} i}}  \right)$ is
obtained from ${\rm {\bf A}}$   by replacing the $j$th column with
the $i$th column, and then by deleting both the $i$th row and
column, $k = \min {\left\{ {I_{n}}  \right.} \setminus {\left.
{\{i\}} \right\}} $.
\end{lemma}
\begin{lemma}\label{lemma:L_ij} \cite{ky3}
Let $L_{i{\kern 1pt} j} $ be the left $ij$-th cofactor of
 ${\rm {\bf A}}\in {\rm M}\left( {n,{\mathbb{H}}} \right)$, that
 is,
$ {\rm{cdet}} _{{j}}\, {\rm {\bf A}} = {{\sum\limits_{i = 1}^{n}
{L_{i{\kern 1pt} j} \cdot a_{i{\kern 1pt} j}} }}$ for all $j
=1,\ldots,n$. Then
\[
L_{i{\kern 1pt} j} = {\left\{ {\begin{array}{*{20}c}
 -{\rm{cdet}} _{i}\, {\rm {\bf A}}_{i{\kern 1pt} .}^{j{\kern 1pt}j} \left( {{\rm {\bf a}}_{j{\kern 1pt}. } }\right),& {i \ne
j},\\
 {\rm{cdet}} _{k}\, {\rm {\bf A}}^{j\, j},& {i = j},
\\
\end{array} }\right.}
\]
\noindent where  ${\rm {\bf A}}_{i{\kern 1pt} .}^{jj} \left( {{\rm
{\bf a}}_{j{\kern 1pt} .} } \right)$ is obtained from ${\rm {\bf
A}}$
 by replacing the $i$th row with the $j$th row, and then by
deleting both the $j$th row and  column, $k = \min {\left\{
{J_{n}} \right.} \setminus {\left. {\{j\}} \right\}} $.
\end{lemma}
The following theorem has a key value in the theory of the column
and row determinants.
\begin{theorem} \cite{ky3}\label{theorem:
determinant of hermitian matrix}
If ${\rm {\bf A}} = \left( {a_{ij}}  \right) \in {\rm M}\left(
{n,{\rm {\mathbb{H}}}} \right)$ is Hermitian, then ${\rm{rdet}}
_{1} {\rm {\bf A}} = \cdots = {\rm{rdet}} _{n} {\rm {\bf A}} =
{\rm{cdet}} _{1} {\rm {\bf A}} = \cdots = {\rm{cdet}} _{n} {\rm
{\bf A}} \in {\rm {\mathbb{R}}}.$
\end{theorem}
\begin{remark}
\label{remark: determinant of hermitian matrix}
 Since all  column and  row determinants of a
Hermitian matrix over ${\rm {\mathbb{H}}}$ are equal, we can
define the determinant of a  Hermitian matrix ${\rm {\bf A}}\in
{\rm M}\left( {n,{\rm {\mathbb{H}}}} \right)$. By definition, we
put for all $i =1,\ldots,n$
\[\det {\rm {\bf A}}: = {\rm{rdet}}_{{i}}\,
{\rm {\bf A}} = {\rm{cdet}} _{{i}}\, {\rm {\bf A}}. \]
\end{remark}
Properties of the determinant of a Hermitian matrix is completely
explored in
 \cite{ky3} by its row and
column determinants. They can be summarized by the following
theorems.
\begin{theorem}\label{theorem:row_combin} If the $i$th row of
a Hermitian matrix ${\rm {\bf A}}\in {\rm M}\left( {n,{\rm
{\mathbb{H}}}} \right)$ is replaced with a left linear combination
of its other rows, i.e. ${\rm {\bf a}}_{i.} = c_{1} {\rm {\bf
a}}_{i_{1} .} + \ldots + c_{k}  {\rm {\bf a}}_{i_{k} .}$, where $
c_{l} \in {{\rm {\mathbb{H}}}}$ for all $ l = 1,\ldots, k$ and
$\{i,i_{l}\}\subset I_{n} $, then
\[
 {\rm{rdet}}_{i}\, {\rm {\bf A}}_{i \, .} \left(
{c_{1} {\rm {\bf a}}_{i_{1} .} + \ldots + c_{k} {\rm {\bf
a}}_{i_{k} .}}  \right) = {\rm{cdet}} _{i}\, {\rm {\bf A}}_{i\, .}
\left( {c_{1}
 {\rm {\bf a}}_{i_{1} .} + \ldots + c_{k} {\rm {\bf
a}}_{i_{k} .}}  \right) = 0.
\]
\end{theorem}
\begin{theorem}\label{theorem:colum_combin} If the $j$th column of
 a Hermitian matrix ${\rm {\bf A}}\in
{\rm M}\left( {n,{\rm {\mathbb{H}}}} \right)$   is replaced with a
right linear combination of its other columns, i.e. ${\rm {\bf
a}}_{.j} = {\rm {\bf a}}_{.j_{1}}   c_{1} + \ldots + {\rm {\bf
a}}_{.j_{k}} c_{k} $, where $c_{l} \in{{\rm {\mathbb{H}}}}$ for
all $l =1,\ldots, k$ and $\{j,j_{l}\}\subset J_{n}$, then
 \[{\rm{cdet}} _{j}\, {\rm {\bf A}}_{.j}
\left( {{\rm {\bf a}}_{.j_{1}} c_{1} + \ldots + {\rm {\bf
a}}_{.j_{k}}c_{k}} \right) ={\rm{rdet}} _{j} \,{\rm {\bf A}}_{.j}
\left( {{\rm {\bf a}}_{.j_{1}}  c_{1} + \ldots + {\rm {\bf
a}}_{.j_{k}}  c_{k}} \right) = 0.
\]
\end{theorem}
The following theorem on the determinantal representation of the
inverse matrix of the Hermitian follows directly from these
properties.
\begin{theorem}\cite{ky3} \label{inver_her} If for a Hermitian matrix ${\rm {\bf A}}\in {\rm M}\left( {n,{\rm
{\mathbb{H}}}} \right)$,
\[\det {\rm {\bf A}} \ne 0,\]
then there exist a unique right
inverse  matrix $(R{\rm {\bf A}})^{ - 1}$ and a unique left
inverse matrix $(L{\rm {\bf A}})^{ - 1}$ of a nonsingular
 ${\rm {\bf A}}$, where $\left( {R{\rm {\bf A}}} \right)^{ - 1} = \left( {L{\rm {\bf A}}}
\right)^{ - 1} = :{\rm {\bf A}}^{ - 1}$, and the right  and left
inverse matrices possess the following determinantal representations

\begin{equation}\label{eq:inver_her_R}
  \left( {R{\rm {\bf A}}} \right)^{ - 1} = {\frac{{1}}{{\det {\rm
{\bf A}}}}}
\begin{pmatrix}
  R_{11} & R_{21} & \cdots & R_{n1}\\
  R_{12} & R_{22} & \cdots & R_{n2}\\
  \cdots & \cdots & \cdots& \cdots\\
  R_{1n} & R_{2n} & \cdots & R_{nn}
\end{pmatrix},
\end{equation}
\begin{equation}\label{eq:inver_her_L}
  \left( {L{\rm {\bf A}}} \right)^{ - 1} = {\frac{{1}}{{\det {\rm
{\bf A}}}}}
\begin{pmatrix}
  L_{11} & L_{21} & \cdots & L_{n1} \\
  L_{12} & L_{22} & \cdots & L_{n2} \\
  \cdots & \cdots & \cdots & \cdots \\
  L_{1n} & L_{2n} & \cdots & L_{nn}
\end{pmatrix},
\end{equation}
 where $R_{ij}$,  $L_{ij}$ are right and left $ij$th cofactors of
${\rm {\bf
 A}}$ respectively for all $ i,j =
1,...,n$.
\end{theorem}

Since the principal submatrices of a Hermitian matrix are Hermitian, the principal minor may be defined as the determinant of its principal submatrix by analogy to the commutative case. We have introduced in \cite{ky4} the rank by principle minors that is the maximal order of a nonzero principal minor of a Hermitian matrix. The following theorem determines a relationship between it and the column
rank of a matrix defining as ceiling amount of right-linearly independent
columns, and the row
rank defining as ceiling amount of left-linearly independent rows.

\begin{theorem} \cite{ky4}\label{theor:rank_her} If  ${\rm {\bf A}}\in {\rm M}\left( {n,{\rm
{\mathbb{H}}}} \right)$ is Hermitian,  then its rank by principal minors are equal to its column and row
ranks.
\end{theorem}

 Due to the noncommutativity of quaternions, there are two
types of eigenvalues.
A quaternion $\lambda$ is said to be a right eigenvalue
of ${\rm {\bf A}} \in {\rm M}\left( {n,{\rm {\mathbb{H}}}}
\right)$ if ${\rm {\bf A}} \cdot {\rm {\bf x}} = {\rm
{\bf x}} \cdot \lambda $ for some nonzero quaternion column-vector
${\rm {\bf x}}$ with quaternion components. Similarly $\lambda$ is a left eigenvalue if ${\rm
{\bf A}} \cdot {\rm {\bf x}} = \lambda \cdot {\rm {\bf x}}$ for some nonzero quaternion column-vector
${\rm {\bf x}}$ with quaternion components.

The theory on the left eigenvalues of quaternion matrices has been
investigated in particular in \cite{hu, so, wo}. The theory on the
right eigenvalues of quaternion matrices is more developed. In
particular we note  \cite{br,ma,ba,dra,zh,far}.

\begin{proposition}\cite{zh}
Let ${\rm {\bf A}} \in {\rm M}\left( {n,{\rm {\mathbb{H}}}}
\right)$ is Hermitian. Then ${\rm {\bf A}}$ has exactly $n$ real
right eigenvalues.
\end{proposition}

Right and left eigenvalues are in general unrelated \cite{fa}, but it is not  for Hermitian matrices.
Suppose ${\rm {\bf A}} \in {\rm M}\left( {n, {\mathbb{H}}}\right)$
is Hermitian and $\lambda \in {\rm {\mathbb {R}}}$ is its right
eigenvalue, then ${\rm {\bf A}} \cdot {\rm {\bf x}} = {\rm {\bf
x}} \cdot \lambda = \lambda \cdot {\rm {\bf x}}$. This means that
all right eigenvalues of a Hermitian matrix are its left
eigenvalues as well. For real left eigenvalues, $\lambda \in {\rm
{\mathbb {R}}}$, the matrix $\lambda {\rm {\bf I}} - {\rm {\bf
A}}$ is Hermitian.
\begin{definition}
If $t \in {\rm {\mathbb {R}}}$, then for a Hermitian matrix ${\rm
{\bf A}}$ the polynomial $p_{{\rm {\bf A}}}\left( {t} \right) =
\det \left( {t{\rm {\bf I}} - {\rm {\bf A}}} \right)$ is said to
be the characteristic polynomial of ${\rm {\bf A}}$.
\end{definition}
The roots of the characteristic polynomial of a Hermitian matrix
are its real left eigenvalues, which are its right eigenvalues as
well. We can prove the following theorem by analogy to the
commutative case (see, e.g. \cite{la}).

\begin{theorem}\label{theor:char_polin}
If ${\rm {\bf A}} \in {\rm M}\left( {n,{\rm {\mathbb{H}}}}
\right)$ is Hermitian, then $p_{{\rm {\bf A}}}\left( {t} \right) =
t^{n} - d_{1} t^{n - 1} + d_{2} t^{n - 2} - \ldots + \left( { - 1}
\right)^{n}d_{n}$, where $d_{k} $ is the sum of principle minors
of ${\rm {\bf A}}$ of order $rk$, $1 \le k < n$, and $d_{n}=\det
{\rm {\bf A}}$.
\end{theorem}

\section{An analogue of the classical adjoint matrix for the Drazin
inverse}
For any  matrix  ${\rm {\bf A}}\in {\mathbb H}^{n\times
n} $ with $ Ind{\kern 1pt} {\rm {\bf A}}=k$, where  a positive
integer $k =: Ind{\kern 1pt} {\rm {\bf A}}= {\mathop {\min}
\limits_{k \in N \cup {\left\{ {0} \right\}}} }{\kern 1pt}
{\left\{ {\rank{\rm {\bf A}}^{k + 1} = \rank{\rm {\bf A}}^{k}}
\right\}}$,   \textbf{the Drazin inverse}  is the unique matrix ${\rm {\bf
X}}$ that satisfies the following three properties
 \begin{equation}\label{eq:Dr_prop}
\begin{array}{l}
  1)\,\,  {\rm {\bf A}}^{k+1}{\rm
{\bf X}}={\rm {\bf
         A}}^{k};\\
  2)\,\,{\rm {\bf X}}{\rm {\bf A}}{\rm {\bf X}}={\rm {\bf X}};\\
  3)\,\, {\rm
{\bf A}}{\rm {\bf X}}={\rm {\bf X}}{\rm {\bf A}}.
\end{array}
\end{equation}
It is denoted by ${\rm {\bf X}}={\rm {\bf A}}^{D}$.

 In particular, when $Ind{\kern 1pt} {\rm {\bf A}}=1$,
then the matrix ${\rm {\bf X}}$ in (\ref{eq:Dr_prop}) is called
\textbf{the group inverse} and is denoted by ${\rm {\bf X}}={\rm {\bf
A}}^{g }$.

If $Ind{\kern 1pt} {\rm {\bf A}}=0$, then ${\rm {\bf A}}$ is
nonsingular, and ${\rm {\bf A}}^{D}\equiv {\rm {\bf A}}^{-1}$.

\begin{remark} Since the equation 3) of (\ref{eq:Dr_prop}), the equation 1)  can be replaced by follows
\[ 1a)\,\,  {\rm {\bf X}}{\rm {\bf A}}^{k+1}={\rm {\bf
         A}}^{k}.\]
\end{remark}
By analogy to the complex case \cite{ca} the following theorem about the limit representation of the
Drazin inverse can be proved.
\begin{theorem} \cite{ca}\label{theor:lim_rep_dr} If ${\rm {\bf A}}
 \in {\mathbb H}^{n\times n}$ with $ Ind{\kern
1pt} {\rm {\bf A}}=k$, then
\[
{\rm {\bf A}}^{D} = {\mathop {\lim} \limits_{\lambda \to 0}}
\left( {\lambda {\rm {\bf I}}_n + {\rm {\bf A}}^{k + 1}} \right)^{
- 1}{\rm {\bf A}}^{k}= {\mathop {\lim} \limits_{\lambda \to 0}}{\rm {\bf A}}^{k}
\left( {\lambda {\rm {\bf I}}_n + {\rm {\bf A}}^{k + 1}} \right)^{
- 1},
\]
where  $\lambda \in {\mathbb R} _{ +}  $, and ${\mathbb R} _{ +} $
is a set of the real positive numbers.
\end{theorem}
 Denote by ${\rm {\bf a}}_{.j}^{(m)} $ and ${\rm {\bf
a}}_{i.}^{(m)} $ the $j$th column  and the $i$th row of  ${\rm
{\bf A}}^{m} $ respectively.
\begin{lemma} \label{lem:rank_col} If ${\rm {\bf A}} \in {\rm M}\left( {n, {\mathbb{H}}}\right)$ with $ Ind{\kern 1pt} {\rm {\bf A}}=k$, then
\begin{equation}\label{eq:rank_col}
 \rank\,\left( {{\rm {\bf A}}^{ k+1} }
\right)_{.\,i} \left( {{\rm {\bf a}}_{.j}^{ (k)} }  \right) \le
\rank\,\left( {{\rm {\bf A}}^{k+1} } \right).
\end{equation}
\end{lemma}
\textit{Proof}
 We can consider the matrix ${\rm {\bf A}}^{
 k+1}$ as ${\rm {\bf A}}^{
 k}{\rm {\bf A}}$.
Let ${\rm {\bf P}}_{i\,s} \left( {-a_{j\,s}}  \right)\in {\mathbb
H}^{n\times n} $, $(s \ne i )$, be a matrix with $-a_{j\,s} $ in
the $(i, s)$ entry, 1 in all diagonal entries, and 0 in others.
The matrix ${\rm {\bf P}}_{i\,s} \left( {-a_{j\,s}}  \right)\in
{\mathbb H}^{n\times n} $, $(s \ne i )$, is a matrix of a
elementary transformation. It follows that
\[
\left( {{\rm {\bf A}}^{ k} {\rm {\bf A}}} \right)_{.\,i} \left(
{{\rm {\bf a}}_{.\,j}^{ (k)} }  \right) \cdot {\prod\limits_{s \ne
i} {{\rm {\bf P}}_{i\,s} \left( {-a_{j\,s}}  \right) = {\mathop
{\left( {{\begin{array}{*{20}c}
 {{\sum\limits_{s \ne j} {a_{1s}^{ (k)}  a_{s1}} } } \hfill & {\ldots}  \hfill
& {a_{1j}^{ (k)} }  \hfill & {\ldots}  \hfill & {{\sum\limits_{s
\ne
j } {a_{1s}^{ (k)}  a_{sn}}}} \hfill \\
 {\ldots}  \hfill & {\ldots}  \hfill & {\ldots}  \hfill & {\ldots}  \hfill &
{\ldots}  \hfill \\
 {{\sum\limits_{s \ne j} {a_{ns}^{(k)}  a_{s1}} } } \hfill & {\ldots}  \hfill
& {a_{nj}^{(k)} }  \hfill & {\ldots}  \hfill & {{\sum\limits_{s
\ne
j } {a_{ns}^{(k)}  a_{sn}}}} \hfill \\
\end{array}} }
\right)}\limits_{i-th}}}}.
\]
 We have the next factorization of the obtained matrix.
\[
{\mathop {\left( {{\begin{array}{*{20}c}
 {{\sum\limits_{s \ne j} {a_{1s}^{ (k)}  a_{s1}} } } \hfill & {\ldots}  \hfill
& {a_{1j}^{ (k)} }  \hfill & {\ldots}  \hfill & {{\sum\limits_{s
\ne
j } {a_{1s}^{ (k)}  a_{sn}} } } \hfill \\
 {\ldots}  \hfill & {\ldots}  \hfill & {\ldots}  \hfill & {\ldots}  \hfill &
{\ldots}  \hfill \\
 {{\sum\limits_{s \ne j} {a_{ns}^{ (k)}  a_{s1}} } } \hfill & {\ldots}  \hfill
& {a_{nj}^{ (k)} }  \hfill & {\ldots}  \hfill & {{\sum\limits_{s
\ne
j } {a_{ns}^{ (k)}  a_{sn}} } } \hfill \\
\end{array}} } \right)}\limits_{i-th}}  =
\]
\[
 = \left( {{\begin{array}{*{20}c}
 {a_{11}^{ (k)} }  \hfill & {a_{12}^{ (k)} }  \hfill & {\ldots}  \hfill &
{a_{1n}^{ (k)} }  \hfill \\
 {a_{21}^{ (k)} }  \hfill & {a_{22}^{ (k)} }  \hfill & {\ldots}  \hfill &
{a_{2n}^{ (k)} }  \hfill \\
 {\ldots}  \hfill & {\ldots}  \hfill & {\ldots}  \hfill & {\ldots}  \hfill
\\
 {a_{n1}^{ (k)} }  \hfill & {a_{n2}^{ (k)} }  \hfill & {\ldots}  \hfill &
{a_{nn}^{ (k)} }  \hfill \\
\end{array}} } \right){\mathop {\left( {{\begin{array}{*{20}c}
 {a_{11}}  \hfill & {\ldots}  \hfill & {0} \hfill & {\ldots}  \hfill &
{a_{1n}}  \hfill \\
 {\ldots}  \hfill & {\ldots}  \hfill & {\ldots}  \hfill & {\ldots}  \hfill &
{\ldots}  \hfill \\
 {0} \hfill & {\ldots}  \hfill & {1} \hfill & {\ldots}  \hfill & {0} \hfill
\\
 {\ldots}  \hfill & {\ldots}  \hfill & {\ldots}  \hfill & {\ldots}  \hfill &
{\ldots}  \hfill \\
 {a_{n1}}  \hfill & {\ldots}  \hfill & {0} \hfill & {\ldots}  \hfill &
{a_{nn}}  \hfill \\
\end{array}} } \right)}\limits_{i-th}} j-th.
\]
 Denote ${\rm {\bf \tilde {A}}}: = {\mathop
{\left( {{\begin{array}{*{20}c}
 {a_{11}}  \hfill & {\ldots}  \hfill & {0} \hfill & {\ldots}  \hfill &
{a_{1n}}  \hfill \\
 {\ldots}  \hfill & {\ldots}  \hfill & {\ldots}  \hfill & {\ldots}  \hfill &
{\ldots}  \hfill \\
 {0} \hfill & {\ldots}  \hfill & {1} \hfill & {\ldots}  \hfill & {0} \hfill
\\
 {\ldots}  \hfill & {\ldots}  \hfill & {\ldots}  \hfill & {\ldots}  \hfill &
{\ldots}  \hfill \\
 {a_{n1}}  \hfill & {\ldots}  \hfill & {0} \hfill & {\ldots}  \hfill &
{a_{nn}}  \hfill \\
\end{array}} } \right)}\limits_{i-th}} j-th$. The matrix
${\rm {\bf \tilde {A}}}$ is obtained from ${\rm {\bf A}}$ by
replacing all entries of the $j$th row  and the $i$th column with
zeroes except for 1 in the $(i, j)$ entry. Since elementary
transformations of a matrix do not change a rank, then $\rank{\rm
{\bf A}}^{ k+1} _{.\,i} \left( {{\rm {\bf a}}_{.j}^{ (k)} }
\right) \le \min {\left\{ {\rank{\rm {\bf A}}^{ k },\rank{\rm {\bf
\tilde {A}}}} \right\}}$. It is obvious that $\rank{\rm {\bf
\tilde {A}}} \ge \rank\,{\rm {\bf A}} \ge \rank{\rm {\bf A}}^{ k}
=\rank{\rm {\bf A}}^{ k+1 }$. From this the inequality (\ref{eq:rank_col})
follows immediately.
$\blacksquare$

The next lemma is  proved similarly.
\begin{lemma} If ${\rm {\bf A}} \in {\rm M}\left( {n, {\mathbb{H}}}\right)$ with $ Ind{\kern 1pt} {\rm {\bf A}}=k$, then $ \rank\left( {{\rm
{\bf A}}^{ k+1} } \right)_{i\,.} \left( {{\rm {\bf a}}_{j\,.}^{
(m)} }  \right) \le \rank\left( {{\rm {\bf A}}^{k+1} } \right).$
\end{lemma}

We shall use the following notations. Let $\alpha : = \left\{
{\alpha _{1} ,\ldots ,\alpha _{k}} \right\} \subseteq {\left\{
{1,\ldots ,m} \right\}}$ and $\beta : = \left\{ {\beta _{1}
,\ldots ,\beta _{k}} \right\} \subseteq {\left\{ {1,\ldots ,n}
\right\}}$ be subsets of the order $1 \le k \le \min {\left\{
{m,n} \right\}}$. By ${\rm {\bf A}}_{\beta} ^{\alpha} $ denote the
submatrix of ${\rm {\bf A}}$ determined by the rows indexed by
$\alpha$ and the columns indexed by $\beta$. Then ${\rm {\bf
A}}{\kern 1pt}_{\alpha} ^{\alpha}$ denotes the principal submatrix
determined by the rows and columns indexed by $\alpha$.
 If ${\rm {\bf A}} \in {\rm
M}\left( {n,{\rm {\mathbb{H}}}} \right)$ is Hermitian, then by
${\left| {{\rm {\bf A}}_{\alpha} ^{\alpha} } \right|}$ denote the
corresponding principal minor of $\det {\rm {\bf A}}$.
 For $1 \leq k\leq n$, the collection of strictly
increasing sequences of $k$ integers chosen from $\left\{
{1,\ldots ,n} \right\}$ is denoted by $\textsl{L}_{ k,
n}: = {\left\{ {\,\alpha :\alpha = \left( {\alpha _{1} ,\ldots
,\alpha _{k}} \right),\,{\kern 1pt} 1 \le \alpha _{1} \le \ldots
\le \alpha _{k} \le n} \right\}}$.  For fixed $i \in \alpha $ and $j \in
\beta $, let $I_{r,\,m} {\left\{ {i} \right\}}: = {\left\{
{\,\alpha :\alpha \in L_{r,m} ,i \in \alpha}  \right\}}{\rm ,}
\quad J_{r,\,n} {\left\{ {j} \right\}}: = {\left\{ {\,\beta :\beta
\in L_{r,n} ,j \in \beta}  \right\}}$.

\begin{lemma}\label{lem:char_cdet}
 If ${\rm {\bf A}} \in {\rm
M}\left( {n,{\rm {\mathbb{H}}}} \right)$ is Hermitian with
$ Ind{\kern 1pt} {\rm {\bf A}}=k$ and $t \in \mathbb{R}$, then

\begin{equation}
\label{eq:char_cdet} {\rm{cdet}} _{i} \left( {t{\rm {\bf I}} +
{\rm {\bf A}}^{k+1} } \right)_{.{\kern 1pt} i} \left( {{\rm {\bf
a}}_{.j}^{(k)} }  \right) = c_{1}^{\left( {ij} \right)} t^{n - 1}
+ c_{2}^{\left( {ij} \right)} t^{n - 2} + \ldots + c_{n}^{\left(
{ij} \right)},
\end{equation}

\noindent where $c_{n}^{\left( {ij} \right)} = {\rm{cdet}} _{i}
\left( {{\rm {\bf A}}^{k+1} } \right)_{.\,i} \left( {{\rm {\bf
a}}_{.\,j}^{ (k)} } \right)$  and $c_{s}^{\left( {ij} \right)} =
{\sum\limits_{\beta \in J_{s,\,n} {\left\{ {i} \right\}}}
{{\rm{cdet}} _{i} \left( {\left( {{\rm {\bf A}}^{k+1} }
\right)_{.\,i} \left( {{\rm {\bf a}}_{.\,j}^{k} } \right)}
\right){\kern 1pt}  _{\beta} ^{\beta} } }$ for all $s = \overline
{1,n - 1} $, $i,j = \overline {1,n}$.
\end{lemma}
\textit{Proof} Denote by ${\rm {\bf b}}_{.{\kern 1pt} {\kern
1pt} i} $ the $i$-th column of the Hermitian matrix
 ${\rm {\bf A}}^{ k+1}  = :\left( {b_{ij}}\right)_{n\times n}$.
Consider the Hermitian matrix $\left( {t{\rm {\bf I}} + {\rm {\bf
A}}^{ k+1} } \right)_{. {\kern 1pt} i} ({\rm {\bf
b}}_{.{\kern 1pt}  i} ) \in {\rm {\mathbb{H}}}^{n\times
n}$. It differs from $\left( {t{\rm {\bf I}} + {\rm {\bf A}}^{
k+1} } \right)$ an entry $b_{ii} $. Taking into account Theorem
\ref{theor:char_polin} we obtain
\begin{equation}
\label{kyr19} \det \left( {t{\rm {\bf I}} + {\rm {\bf A}}^{ k+1} }
\right)_{.{\kern 1pt} i} \left( {{\rm {\bf b}}_{.{\kern 1pt}
{\kern 1pt} i}}  \right) = d_{1} t^{n - 1} + d_{2} t^{n - 2} +
\ldots + d_{n},
\end{equation}
where $d_{s} = {\sum\limits_{\beta \in J_{s,\,n} {\left\{ {i}
\right\}}} {| \left( {{\rm {\bf A}}^{ k+1} } \right) _{\beta} ^{\beta} }| } $ is the sum of all principal
minors of order $s$ that contain the $i$-th column for all $s =
\overline {1,n - 1} $ and $d_{n} = \det \left( {{\rm {\bf A}}^{
k+1} } \right)$. Consequently we have ${\rm {\bf b}}_{.
{\kern 1pt} i} = \left( {{\begin{array}{*{20}c}
 {{\sum\limits_{l} {a_{1l}^{(k)}  a_{li}} } } \hfill \\
 {{\sum\limits_{l} {a_{2l}^{(k)}  a_{li}} } } \hfill \\
 { \vdots}  \hfill \\
 {{\sum\limits_{l} {a_{nl}^{(k)}  a_{li}} } } \hfill \\
\end{array}} } \right) = {\sum\limits_{l} {{\rm {\bf a}}_{.\,l}^{(k)}  a_{li}}
}$, where ${\rm {\bf a}}_{.{\kern 1pt}  l}^{(k)}  $ is
the $l$th column-vector  of ${\rm {\bf A}}^{k}$ for all
$l=\overline{1,n}$. Taking into account Theorem \ref{theorem:
determinant of hermitian matrix}, Lemma \ref{lemma:L_ij} and
Proposition \ref{pr:col_mult} we obtain on the one hand
\begin{equation}
\label{kyr20}
\begin{array}{c}
  \det \left( {t{\rm {\bf I}} + {\rm {\bf A}}^{ k+1} }
\right)_{.{\kern 1pt} i} \left( {{\rm {\bf b}}_{.{\kern 1pt}
{\kern 1pt} i} } \right) = {\rm{cdet}} _{i} \left( {t{\rm {\bf I}}
+ {\rm {\bf A}}^{ k+1} } \right)_{.{\kern 1pt} i}
\left( {{\rm {\bf b}}_{.{\kern 1pt} {\kern 1pt} i}}  \right) =\\
   = {\sum\limits_{l} {{\rm{cdet}} _{i} \left( {t{\rm {\bf I}} + {\rm {\bf
   A}}^{k+1
}} \right)_{.{\kern 1pt} l} \left( {{\rm {\bf a}}_{. {\kern 1pt}
l}^{(k)}  a_{l{\kern 1pt} i}} \right) = {\sum\limits_{l}
{{\rm{cdet}} _{i} \left( {t{\rm {\bf I}} + {\rm {\bf A}}^{ k+1} }
\right)_{.{\kern 1pt} i} \left( {{\rm {\bf a}}_{.{\kern 1pt}
{\kern 1pt} l}^{(k)} }  \right) \cdot {\kern 1pt}} } }} a_{li}
\end{array}
\end{equation}
On the other hand having changed the order of summation,  we get
for all $s = \overline {1,n - 1} $
\begin{equation}
\label{kyr21}
\begin{array}{c}
  d_{s} = {\sum\limits_{\beta \in J_{s,\,n} {\left\{ {i} \right\}}}
{\det \left( {{\rm {\bf A}}^{ k+1} } \right){\kern 1pt} {\kern
1pt} _{\beta} ^{\beta} } }  = {\sum\limits_{\beta \in J_{s,\,n}
{\left\{ {i} \right\}}} {{\rm{cdet}} _{i} \left( {{\rm {\bf A}}^{
k+1} } \right){\kern 1pt} {\kern 1pt}
_{\beta} ^{\beta} } }  =\\
  {\sum\limits_{\beta \in J_{s,\,n}
{\left\{ {i} \right\}}} {{\sum\limits_{l} {{\rm{cdet}} _{i} \left(
{\left( {{\rm {\bf A}}^{ k+1} {\kern 1pt}} \right)_{.\, i} \left(
{{\rm {\bf a}}_{.\, l}^{(k)} a_{l\,i}} \right)} \right)}} }} \,
_{\beta} ^{\beta}=\\
  {\sum\limits_{l} { {{\sum\limits_{\beta \in J_{s,\,n} {\left\{ {i}
\right\}}} {{\rm{cdet}} _{i} \left( {\left( {{\rm {\bf A}}^{ k+1}
{\kern 1pt}}  \right)_{.{\kern 1pt} i} \left( {{\rm {\bf
a}}_{.{\kern 1pt} {\kern 1pt} l}^{(k)} }  \right)} \right){\kern
1pt} _{\beta} ^{\beta} } }} }}  \cdot a_{l{\kern 1pt} i}.
\end{array}
\end{equation}
By substituting (\ref{kyr20}) and (\ref{kyr21}) in (\ref{kyr19}),
and equating factors at $a _ {l \, i} $ when $l = j $, we obtain
the equality (\ref{eq:char_cdet}).
$\blacksquare$

 By analogy can be proved
the following lemma.
\begin{lemma}
If ${\rm {\bf A}} \in {\rm M}\left( {n, {\mathbb{H}}}\right)$  is Hermitian with
$ Ind{\kern 1pt} {\rm {\bf A}}=k$ and $t \in \mathbb{R}$, then
\[ {\rm{rdet}} _{j}  {( t{\rm {\bf I}} + {\rm {\bf A}}^{k+1} )_{j\,.\,} ({\rm {\bf a}}_{i.}^{(k)}
)}
   = r_{1}^{\left( {ij} \right)} t^{n
- 1} +r_{2}^{\left( {ij} \right)} t^{n - 2} + \ldots +
r_{n}^{\left( {ij} \right)},
\]

\noindent where  $r_{n}^{\left( {ij} \right)} = {\rm{rdet}} _{j}
{({\rm {\bf A}}^{k+1} )_{j\,.\,} ({\rm {\bf a}}_{i.\,}^{(k)} )}$
and $r_{s}^{\left( {ij} \right)} = {{{\sum\limits_{\alpha \in
I_{s,n} {\left\{ {j} \right\}}} {{\rm{rdet}} _{j} \left( {({\rm
{\bf A}}^{k+1} )_{j\,.\,} ({\rm {\bf a}}_{i.\,}^{(k)} )}
\right)\,_{\alpha} ^{\alpha} } }}}$ for all $s = \overline {1,n -
1} $ and $i,j = \overline {1,n}$.
\end{lemma}

\begin{theorem}\label{theor:det_rep_draz}
If ${\rm {\bf A}} \in {\rm M}\left( {n, {\mathbb{H}}}\right)$ is Hermitian with
$ Ind{\kern 1pt} {\rm {\bf A}}=k$ and $\rank{\rm {\bf A}}^{k+1} =
\rank{\rm {\bf A}}^{k} = r$, then the Drazin inverse ${\rm {\bf
A}}^{D} = \left( {a_{ij}^{D} } \right) \in {\rm
{\mathbb{H}}}_{}^{n\times n} $ possess the following determinantal
representations:
\begin{equation}
\label{eq:dr_rep_cdet} a_{ij}^{D}  = {\frac{{{\sum\limits_{\beta
\in J_{r,\,n} {\left\{ {i} \right\}}} {{\rm{cdet}} _{i} \left(
{\left( {{\rm {\bf A}}^{k+1}} \right)_{\,. \,i} \left( {{\rm {\bf
a}}_{.j}^{ k} }  \right)} \right){\kern 1pt} {\kern 1pt} _{\beta}
^{\beta} } } }}{{{\sum\limits_{\beta \in J_{r,\,\,n}} {{\left|
{\left( {{\rm {\bf A}}^{k+1}} \right){\kern 1pt} _{\beta} ^{\beta}
}  \right|}}} }}},
\end{equation}
or
\begin{equation}
\label{eq:dr_rep_rdet} a_{ij}^{D}  = {\frac{{{\sum\limits_{\alpha
\in I_{r,n} {\left\{ {j} \right\}}} {{\rm{rdet}} _{j} \left(
{({\rm {\bf A}}^{ k+1} )_{j\,.\,} ({\rm {\bf a}}_{i.\,}^{ (k)} )}
\right)\,_{\alpha} ^{\alpha} } }}}{{{\sum\limits_{\alpha \in
I_{r,\,n}}  {{\left| {\left( {{\rm {\bf A}}^{k+1} } \right){\kern
1pt}  _{\alpha} ^{\alpha} } \right|}}} }}}.
\end{equation}
\end{theorem}
\textit{Proof}   At first we prove (\ref{eq:dr_rep_cdet}). By
Theorem \ref{theor:lim_rep_dr}, ${\rm {\bf A}}^{ +}  = {\mathop
{\lim} \limits_{\alpha \to 0}} \left( {\alpha {\rm {\bf I}}_{n} +
{\rm {\bf A}}^{k+1}} \right)^{ - 1}{\rm {\bf A}}^{ k }$. The
matrix $\left( {\alpha {\rm {\bf I}} + {\rm {\bf A}}^{k+1}}
\right) \in {\rm {\mathbb{H}}}^{n\times n}$ is a full-rank
Hermitian matrix. Taking into account Theorem \ref{inver_her} it
has an inverse, which we represent as a left inverse matrix

\[
\left( {\alpha {\rm {\bf I}} + {\rm {\bf A}}^{k+1}} \right)^{ - 1}
= {\frac{{1}}{{\det \left( {\alpha {\rm {\bf I}} + {\rm {\bf A}}^{
k+1}} \right)}}}\left( {{\begin{array}{*{20}c}
 {L_{11}}  \hfill & {L_{21}}  \hfill & {\ldots}  \hfill & {L_{n1}}  \hfill
\\
 {L_{12}}  \hfill & {L_{22}}  \hfill & {\ldots}  \hfill & {L_{n2}}  \hfill
\\
 {\ldots}  \hfill & {\ldots}  \hfill & {\ldots}  \hfill & {\ldots}  \hfill
\\
 {L_{1n}}  \hfill & {L_{2n}}  \hfill & {\ldots}  \hfill & {L_{nn}}  \hfill
\\
\end{array}} } \right),
\]
\noindent where $L_{ij} $ is a left $ij$-th cofactor  of a matrix
$\alpha {\rm {\bf I}} + {\rm {\bf A}}^{k+1}$. Then we have
\[\begin{array}{l}
   \left( {\alpha {\rm {\bf I}} + {\rm {\bf A}}^{k+1}} \right)^{ -
1}{\rm {\bf A}}^{ k}  = \\
  ={\frac{{1}}{{\det \left( {\alpha {\rm
{\bf I}} + {\rm {\bf A}}^{m+1}} \right)}}}\left(
{{\begin{array}{*{20}c}
 {{\sum\limits_{s = 1}^{n} {L_{s1} a_{s1}^{ (k)} } } } \hfill &
{{\sum\limits_{s = 1}^{n} {L_{s1} a_{s2}^{ (k)} } } } \hfill &
{\ldots}
\hfill & {{\sum\limits_{s = 1}^{n} {L_{s1} a_{sn}^{ (k)} } } } \hfill \\
 {{\sum\limits_{s = 1}^{n} {L_{s2} a_{s1}^{ (k)} } } } \hfill &
{{\sum\limits_{s = 1}^{n} {L_{s2} a_{s2}^{ (k)} } } } \hfill &
{\ldots}
\hfill & {{\sum\limits_{s = 1}^{n} {L_{s2} a_{sn}^{ (k)} } } } \hfill \\
 {\ldots}  \hfill & {\ldots}  \hfill & {\ldots}  \hfill & {\ldots}  \hfill
\\
 {{\sum\limits_{s = 1}^{n} {L_{sn} a_{s1}^{ (k)} } } } \hfill &
{{\sum\limits_{s = 1}^{n} {L_{sn} a_{s2}^{ (k)} } } } \hfill &
{\ldots}
\hfill & {{\sum\limits_{s = 1}^{n} {L_{sn} a_{sn}^{ (k)} } } } \hfill \\
\end{array}} } \right).
\end{array}
\]
By using the definition of a left cofactor, we obtain
\begin{equation}
\label{eq:A_d1} {\rm {\bf A}}^{D}  = {\mathop {\lim}
\limits_{\alpha \to 0}} \left( {{\begin{array}{*{20}c}
 {{\frac{{{\rm cdet} _{1} \left( {\alpha {\rm {\bf I}} + {\rm {\bf A}}^{k+1}} \right)_{.1} \left( {{\rm {\bf a}}_{.1}^{ (k)} }
\right)}}{{\det \left( {\alpha {\rm {\bf I}} + {\rm {\bf
A}}^{k+1}} \right)}}}} \hfill & {\ldots}  \hfill & {{\frac{{{\rm
cdet} _{1} \left( {\alpha {\rm {\bf I}} + {\rm {\bf A}}^{k+1}}
\right)_{.1} \left( {{\rm {\bf a}}_{.n}^{(k)} } \right)}}{{\det
\left( {\alpha {\rm {\bf I}} + {\rm {\bf
A}}^{k+1}} \right)}}}} \hfill \\
 {\ldots}  \hfill & {\ldots}  \hfill & {\ldots}  \hfill \\
 {{\frac{{{\rm cdet} _{n} \left( {\alpha {\rm {\bf I}} + {\rm {\bf A}}^{k+1}} \right)_{.n} \left( {{\rm {\bf a}}_{.1}^{ (k)} }
\right)}}{{\det \left( {\alpha {\rm {\bf I}} + {\rm {\bf
A}}^{k+1}} \right)}}}} \hfill & {\ldots}  \hfill & {{\frac{{{\rm
cdet} _{n} \left( {\alpha {\rm {\bf I}} + {\rm {\bf A}}^{k+1}}
\right)_{.n} \left( {{\rm {\bf a}}_{.n}^{(k)} } \right)}}{{\det
\left( {\alpha {\rm {\bf I}} + {\rm {\bf
A}}^{k+1}} \right)}}}} \hfill \\
\end{array}} } \right).
\end{equation}
By Theorem \ref{theor:char_polin} we have
 \[\det \left( {\alpha
{\rm {\bf I}} + {\rm {\bf A}}^{m+1}} \right) = \alpha ^{n} + d_{1}
\alpha ^{n - 1} + d_{2} \alpha ^{n - 2} + \ldots + d_{n},\]
where
$d_{s} = {\sum\limits_{\beta \in J_{s,\,n}}  {{\left| {\left(
{{\rm {\bf A}}^{k+1}} \right){\kern 1pt} {\kern 1pt} _{\beta}
^{\beta} } \right|}}} $ is a sum of principal minors of ${\rm {\bf
A}}^{k+1}$ of order $s$ for all  $s = \overline {1,n - 1} $ and
$d_{n} = \det {\rm {\bf A}}^{k+1}$.

Since $\rank{\rm {\bf A}}^{k+1} = \rank{\rm {\bf A}}^{k} = r$,
then $d_{n} = d_{n - 1} = \ldots = d_{r + 1} = 0$. It follows that
$\det \left( {\alpha {\rm {\bf I}} + {\rm {\bf A}}^{k+1}} \right)
= \alpha ^{n} + d_{1} \alpha ^{n - 1} + d_{2} \alpha ^{n - 2} +
\ldots + d_{r} \alpha ^{n - r}$.

Using (\ref{eq:char_cdet})  we have
 \[{\rm{cdet}} _{i} \left(
{\alpha {\rm {\bf I}} +{\rm {\bf A}}^{k+1} }
\right)_{.i} \left( {{\rm {\bf a}}_{.j}^{ (k)} } \right) =
c_{1}^{\left( {ij} \right)} \alpha ^{n - 1} + c_{2}^{\left( {ij}
\right)} \alpha ^{n - 2} + \ldots + c_{n}^{\left( {ij} \right)} \]
for all $i,j = \overline {1,n} $, where $c_{s}^{\left( {ij}
\right)} = {\sum\limits_{\beta \in J_{s,\,n} {\left\{ {i}
\right\}}} {{\rm{cdet}} _{i} \left( {({\rm {\bf A}}^{k+1})_{.\,i}
\left( {{\rm {\bf a}}_{.j}^{ (k)} } \right)} \right){\kern 1pt}
{\kern 1pt} _{\beta }^{\beta} } } $ for all $s = \overline {1,n -
1} $ and $c_{n}^{\left( {ij} \right)} = {\rm{cdet}} _{i} \left(
{{\rm {\bf A}}^{k+1}} \right)_{.i} \left( {{\rm {\bf a}}_{.j}^{ (k)}
} \right)$.

We prove that $c_{k}^{\left( {ij} \right)} = 0$, when $k \ge r +
1$ for all $i,j = \overline {1,n} $. By Lemma \ref{lem:rank_col}
$\left( {{\rm {\bf A}}^{k+1}} \right)_{.\,i}
\left( {{\rm {\bf a}}_{.j}^{ (k)} }  \right) \le r$, then the matrix
$\left( {{\rm {\bf A}}^{k+1}} \right)_{.\,i} \left( {{\rm {\bf
a}}_{.j}^{ (k)} }  \right)$ has no more $r$ right-linearly
independent columns.

Consider $\left( {({\rm {\bf A}}^{k+1})_{\,.\,i} \left( {{\rm {\bf
a}}_{.j}^{ (k)} } \right)} \right){\kern 1pt} {\kern 1pt} _{\beta}
^{\beta}  $, when $\beta \in J_{s,n} {\left\{ {i} \right\}}$. It
is a principal submatrix of  $\left( {{\rm {\bf A}}^{k+1}}
\right)_{.\,i} \left( {{\rm {\bf a}}_{.j}^{ (k)} } \right)$ of order
$s \ge r + 1$. Deleting both its $i$-th row and column, we obtain
a principal submatrix of order $s - 1$ of  ${\rm {\bf A}}^{k+1}$.
 We denote it by ${\rm {\bf M}}$. The following cases are
possible.

Let $s = r + 1$ and $\det {\rm {\bf M}} \ne 0$. In this case all
columns of $ {\rm {\bf M}} $ are right-linearly independent. The
addition of all of them on one coordinate to columns of
 $\left( {\left( {{\rm {\bf A}}^{k+1}}
\right)_{.\,i} \left( {{\rm {\bf a}}_{.j}^{ (k)} }  \right)} \right)
{\kern 1pt} _{\beta} ^{\beta}$ keeps their right-linear
independence. Hence, they are basis in a matrix $\left( {\left(
{{\rm {\bf A}}^{k+1}} \right)_{\,.\,i} \left( {{\rm {\bf
a}}_{.j}^{ (k)} } \right)} \right){\kern 1pt} {\kern 1pt} _{\beta}
^{\beta} $, and  the $i$-th
column is the right linear combination of its basis columns. From
this by Theorem \ref{theorem:colum_combin}, we get ${\rm{cdet}}
_{i} \left( {\left( {{\rm {\bf A}}^{k+1}} \right)_{\,.\,i} \left(
{{\rm {\bf a}}_{.j}^{ (k)} } \right)} \right){\kern 1pt} {\kern 1pt}
_{\beta} ^{\beta}  = 0$, when $\beta \in J_{s,n} {\left\{ {i}
\right\}}$ and $s= r + 1$.

If $s = r + 1$ and $\det {\rm {\bf M}} = 0$, than $p$, ($p < s$),
columns are basis in  ${\rm {\bf M}}$ and in  $\left( {\left(
{{\rm {\bf A}}^{k+1}} \right)_{.\,i} \left( {{\rm {\bf a}}_{.j}^{
(k)} } \right)} \right){\kern 1pt} {\kern 1pt} _{\beta} ^{\beta} $.
Then by Theorems \ref{theor:rank_her} and \ref{theorem:colum_combin}  we obtain ${\rm{cdet}} _{i} \left(
{\left( {{\rm {\bf A}}^{k+1}} \right)_{\,.\,i} \left( {{\rm {\bf
a}}_{.j}^{ (k)} } \right)} \right){\kern 1pt} {\kern 1pt} _{\beta}
^{\beta}  = 0$ as well.

If $s > r + 1$, then from Theorem \ref{theor:rank_her}  it follows that $\det {\rm {\bf M}} = 0$ and
$p$, ($p < r$), columns are basis in the both matrices ${\rm {\bf
M}}$ and $\left( {\left( {{\rm {\bf A}}^{k+1}} \right)_{\,.\,i}
\left( {{\rm {\bf a}}_{.j}^{ (k)} }  \right)} \right){\kern 1pt}
{\kern 1pt} _{\beta} ^{\beta}  $. Then by Theorem
\ref{theorem:colum_combin}, we
have ${\rm{cdet}} _{i} \left( {\left( {{\rm {\bf A}}^{k+1}}
\right)_{\,.\,i} \left( {{\rm {\bf a}}_{.j}^{ (k)} } \right)}
\right){\kern 1pt} {\kern 1pt} _{\beta} ^{\beta}  = 0.$

Thus in all cases we have ${\rm{cdet}} _{i} \left( {\left( {{\rm
{\bf A}}^{k+1}} \right)_{\,.\,i} \left( {{\rm {\bf a}}_{.j}^{ (k)} }
\right)} \right){\kern 1pt} {\kern 1pt} _{\beta} ^{\beta}  = 0$,
when $\beta \in J_{s,n} {\left\{ {i} \right\}}$ and $r + 1 \le s <
n$. From here if $r + 1 \le s < n$, then
\[c_{s}^{\left( {ij} \right)} = {\sum\limits_{\beta \in J_{s,\,n}
{\left\{ {i} \right\}}} {{\rm{cdet}} _{i} \left( {\left( {{\rm
{\bf A}}^{k+1}} \right)_{\,.\,i} \left( {{\rm {\bf a}}_{.j}^{ (k)} }
\right)} \right) {\kern 1pt} _{\beta} ^{\beta} } }=0,\] and
$c_{n}^{\left( {ij} \right)} = {\rm{cdet}} _{i} \left( {{\rm {\bf
A}}^{k+1}} \right)_{.\,i} \left( {{\rm {\bf a}}_{.j}^{(k)} }
\right) = 0$ for all $i,j = \overline {1,n} $.

Hence, ${\rm{cdet}} _{i} \left( {\alpha {\rm {\bf I}} + {\rm {\bf
A}}^{k+1}} \right)_{.\,i} \left( {{\rm {\bf a}}_{.\,j}^{ (k)} }
\right) =c_{1}^{\left( {ij} \right)} \alpha ^{n - 1} +
c_{2}^{\left( {ij} \right)} \alpha ^{n - 2} + \ldots +
c_{r}^{\left( {ij} \right)} \alpha ^{n - r}$ for all $i,j =
\overline {1,n} $. By substituting these values in the matrix from
(\ref{eq:A_d1}), we obtain
\[\begin{array}{c}
  {\rm {\bf A}}^{ D}  = {\mathop {\lim} \limits_{\alpha \to 0}} \left(
{{\begin{array}{*{20}c}
 {{\frac{{c_{1}^{\left( {11} \right)} \alpha ^{n - 1} + \ldots +
c_{r}^{\left( {11} \right)} \alpha ^{n - r}}}{{\alpha ^{n} + d_{1}
\alpha ^{n - 1} + \ldots + d_{r} \alpha ^{n - r}}}}} \hfill &
{\ldots}  \hfill & {{\frac{{c_{1}^{\left( {1n} \right)} \alpha ^{n
- 1} + \ldots + c_{r}^{\left( {1n} \right)} \alpha ^{n -
r}}}{{\alpha ^{n} + d_{1} \alpha
^{n - 1} + \ldots + d_{r} \alpha ^{n - r}}}}} \hfill \\
 {\ldots}  \hfill & {\ldots}  \hfill & {\ldots}  \hfill \\
 {{\frac{{c_{1}^{\left( {n1} \right)} \alpha ^{n - 1} + \ldots +
c_{r}^{\left( {n1} \right)} \alpha ^{n - r}}}{{\alpha ^{n} + d_{1}
\alpha ^{n - 1} + \ldots + d_{r} \alpha ^{n - r}}}}} \hfill &
{\ldots}  \hfill & {{\frac{{c_{1}^{\left( {nn} \right)} \alpha ^{n
- 1} + \ldots + c_{r}^{\left( {nn} \right)} \alpha ^{n -
r}}}{{\alpha ^{n} + d_{1} \alpha
^{n - 1} + \ldots + d_{r} \alpha ^{n - r}}}}} \hfill \\
\end{array}} } \right) =\\
  \left( {{\begin{array}{*{20}c}
 {{\frac{{c_{r}^{\left( {11} \right)}} }{{d_{r}} }}} \hfill & {\ldots}
\hfill & {{\frac{{c_{r}^{\left( {1n} \right)}} }{{d_{r}} }}} \hfill \\
 {\ldots}  \hfill & {\ldots}  \hfill & {\ldots}  \hfill \\
 {{\frac{{c_{r}^{\left( {n1} \right)}} }{{d_{r}} }}} \hfill & {\ldots}
\hfill & {{\frac{{c_{r}^{\left( {nn} \right)}} }{{d_{r}} }}} \hfill \\
\end{array}} } \right).
\end{array}
\]

Here $c_{r}^{\left( {ij} \right)} = {\sum\limits_{\beta \in
J_{r,\,n} {\left\{ {i} \right\}}} {{\rm{cdet}} _{i} \left( {\left(
{{\rm {\bf A}}^{k+1}} \right)_{\,.\,i} \left( {{\rm {\bf
a}}_{.j}^{ (k)} }  \right)} \right){\kern 1pt} {\kern 1pt} _{\beta}
^{\beta} } } $ and $d_{r} = {\sum\limits_{\beta \in J_{r,\,\,n}}
{{\left| {\left(
{{\rm {\bf A}}^{k+1}} \right){\kern
1pt} {\kern 1pt} _{\beta} ^{\beta} } \right|}}} $. Thus, we have
obtained the determinantal representation of ${\rm {\bf A}}_{}^{
+}  $ by (\ref{eq:dr_rep_cdet}).

By analogy can be proved the determinantal representation of ${\rm
{\bf A}}^{ D} $ by (\ref{eq:dr_rep_rdet}).
$\blacksquare$

 In
the following corollaries we introduce  determinantal
 representations of the group inverse ${\rm {\bf A}}^{g }$ and
  the matrix ${\rm {\bf A}}^{D}{\rm {\bf A}}$ respectively.
\begin{Corollary}
 If $Ind{\kern 1pt} {\rm {\bf A}} = 1 $ and $\rank{\rm {\bf
A}}^{2} = \rank{\rm {\bf A}}=r \le n$ for a Hermitian matrix
${\rm {\bf A}}\in {\mathbb H}^{n\times n} $, then the group
inverse ${\rm {\bf A}}^{g }$ possess the following
determinantal representations:
\[ a_{ij}^{g}  =
{\frac{{{\sum\limits_{\beta \in J_{r,\,n} {\left\{ {i} \right\}}}
{{\rm{cdet}} _{i} \left( {\left( {{\rm {\bf A}}^{2}} \right)_{\,.
\,i} \left( {{\rm {\bf a}}_{.j} } \right)} \right){\kern 1pt}
{\kern 1pt} _{\beta} ^{\beta} } } }}{{{\sum\limits_{\beta \in
J_{r,\,\,n}} {{\left| {\left( {{\rm {\bf A}}^{2}} \right){\kern
1pt} _{\beta} ^{\beta} }  \right|}}} }}},
\]
or
\[ a_{ij}^{g}  = {\frac{{{\sum\limits_{\alpha \in I_{r,n}
{\left\{ {j} \right\}}} {{\rm{rdet}} _{j} \left( {({\rm {\bf A}}^{
2} )_{j\,.\,} ({\rm {\bf a}}_{i.\,} )} \right)\,_{\alpha}
^{\alpha} } }}}{{{\sum\limits_{\alpha \in I_{r,\,n}}  {{\left|
{\left( {{\rm {\bf A}}^{2} } \right){\kern 1pt}  _{\alpha}
^{\alpha} } \right|}}} }}}.
\]
\end{Corollary}
\textit{Proof} The proof follows immediately from Theorem
\ref{theor:det_rep_draz} in view of $k=1$.$\blacksquare$

\begin{Corollary}
If $Ind{\kern 1pt} {\rm {\bf A}} = k $ and $\rank{\rm {\bf A}}^{k
+ 1} = \rank{\rm {\bf A}}^{k}=r \le n$ for an arbitrary matrix
${\rm {\bf A}}\in {\mathbb C}^{n\times n} $, then
\begin{equation}
\label{eq:rep_A_d_A}
  {\rm {\bf A}}^{D}{\rm {\bf A}} = \left({{\frac{{{\sum\limits_{\beta
\in J_{r,\,n} {\left\{ {i} \right\}}} {{\rm{cdet}} _{i} \left(
{\left( {{\rm {\bf A}}^{k+1}} \right)_{\,. \,i} \left( {{\rm {\bf
a}}_{.j}^{ (k+1)} }  \right)} \right){\kern 1pt} {\kern 1pt}
_{\beta} ^{\beta} } } }}{{{\sum\limits_{\beta \in J_{r,\,\,n}}
{{\left| {\left( {{\rm {\bf A}}^{k+1}} \right){\kern 1pt} _{\beta}
^{\beta} }  \right|}}} }}}}\right)_{n\times n},
\end{equation}
and
\begin{equation}
\label{eq:rep_A_A_d}
  {\rm {\bf A}}{\rm {\bf A}}^{D} = \left({{\frac{{{\sum\limits_{\alpha
\in I_{r,n} {\left\{ {j} \right\}}} {{\rm{rdet}} _{j} \left(
{({\rm {\bf A}}^{ k+1} )_{j\,.\,} ({\rm {\bf a}}_{i.\,}^{ (k+1)}
)} \right)\,_{\alpha} ^{\alpha} } }}}{{{\sum\limits_{\alpha \in
I_{r,\,n}}  {{\left| {\left( {{\rm {\bf A}}^{k+1} } \right){\kern
1pt}  _{\alpha} ^{\alpha} } \right|}}} }}}}\right)_{n\times n}.
\end{equation}
\end{Corollary}
\textit{Proof} At first we prove (\ref{eq:rep_A_d_A}). Let
${\rm {\bf A}}^{D}{\rm {\bf A}}=(v_{ij})_{n\times n}$. Using
(\ref{eq:dr_rep_cdet}) for arbitrary $1\leq i,j\leq n$ we have
 \[
  v_{i\,j}=\sum\limits_{s}
 {{\frac{{{\sum\limits_{\beta
\in J_{r,\,n} {\left\{ {i} \right\}}} {{\rm{cdet}} _{i} \left(
{\left( {{\rm {\bf A}}^{k+1}} \right)_{\,. \,i} \left( {{\rm {\bf
a}}_{.j}^{ (k)} }  \right)} \right){\kern 1pt} {\kern 1pt} _{\beta}
^{\beta} } } }}{{{\sum\limits_{\beta \in J_{r,\,\,n}} {{\left|
{\left( {{\rm {\bf A}}^{k+1}} \right){\kern 1pt} _{\beta} ^{\beta}
}  \right|}}} }}}}\cdot a_{s\,j}=\]
  \[{{\frac{{{\sum\limits_{\beta
\in J_{r,\,n} {\left\{ {i} \right\}}}\sum\limits_{s} {{\rm{cdet}}
_{i} \left( {\left( {{\rm {\bf A}}^{k+1}} \right)_{\,. \,i} \left(
{{\rm {\bf a}}_{.j}^{ (k)}\cdot a_{s\,j} }  \right)} \right){\kern
1pt} {\kern 1pt} _{\beta} ^{\beta} } } }}{{{\sum\limits_{\beta \in
J_{r,\,\,n}} {{\left| {\left( {{\rm {\bf A}}^{k+1}} \right){\kern
1pt} _{\beta} ^{\beta} }  \right|}}}
}}}}={{\frac{{{\sum\limits_{\beta \in J_{r,\,n} {\left\{ {i}
\right\}}} {{\rm{cdet}} _{i} \left( {\left( {{\rm {\bf A}}^{k+1}}
\right)_{\,. \,i} \left( {{\rm {\bf a}}_{.j}^{ (k+1)} }  \right)}
\right){\kern 1pt} {\kern 1pt} _{\beta} ^{\beta} } }
}}{{{\sum\limits_{\beta \in J_{r,\,\,n}} {{\left| {\left( {{\rm
{\bf A}}^{k+1}} \right){\kern 1pt} _{\beta} ^{\beta} } \right|}}}
}}}}.
\]
By analogy can be proved (\ref{eq:rep_A_A_d}) using the
determinantal representation of the Drazin inverse by
(\ref{eq:dr_rep_rdet}).
$\blacksquare$

\section{Cramer's rule of the Drazin inverse solutions of some  matrix
equations}
Consider a  matrix equation
\begin{equation}\label{eq:AX=B}
 {\rm {\bf A}}{\rm {\bf X}} = {\rm {\bf B}},
\end{equation}
where $ {\bf
A}\in {\mathbb{H}}^{n\times n} $, $ {\rm {\bf B}}\in {\mathbb{H}}^{n\times m} $ are given, ${\rm {\bf
A}}$ is Hermitian and ${\rm {\bf X}} \in
{\mathbb{H}}^{n\times m}$ is unknown. Let $Ind{\kern 1pt} {\rm {\bf A}} = k $.
 We denote ${\rm {\bf A}}^{ k}{\rm {\bf B}}=:\hat{{\rm
{\bf B}}}= (\hat{b}_{ij})\in {\mathbb{H}}^{n\times m}$.
\begin{theorem}
If $\rank{\rm {\bf
A}}^{k + 1} = \rank{\rm {\bf A}}^{k}=r \le n$ for ${\rm {\bf A}}\in {\mathbb H}^{n\times n} $, then for Drazin inverse solution $ {\bf X}={\bf A}^{D}{\bf B}= (x_{ij})$ of  (\ref{eq:AX=B}) we have
\begin{equation}
\label{eq:dr_AX} x_{ij} = {\frac{{{\sum\limits_{\beta \in
J_{r,\,n} {\left\{ {i} \right\}}} {{{\rm{cdet}} _{i} \left( {\left( {{\rm
{\bf A}}^{ k+1}} \right)_{\,.\,i} \left( {\hat{{\rm
{\bf b}}}_{.j}} \right)} \right){\kern 1pt} {\kern 1pt} _{\beta}
^{\beta} }} } }}{{{\sum\limits_{\beta \in J_{r,\,\,n}}
{{\left| {\left( {{\rm {\bf A}}^{ k+1} } \right){\kern
1pt} {\kern 1pt} _{\beta} ^{\beta} }  \right|}}} }}}.
\end{equation}
\end{theorem}

\textit{Proof}  By Theorem   \ref{theor:det_rep_draz} we can represent the matrix
${\rm {\bf A}}^{ D} $ by ( \ref{eq:dr_rep_cdet}). Therefore, we
obtain for all $i = \overline {1,n} $ and $j = \overline {1,m} $

\[
x_{ij}  =\sum_{s=1}^{n}
a_{is}^{D}b_{sj}=\sum_{s=1}^{n}{\frac{{{\sum\limits_{\beta \in
J_{r,\,n} {\left\{ {i} \right\}}} {{\rm{cdet}} _{i} \left( {\left(
{{\rm {\bf A}}^{k+1} } \right)_{\,. \,i} \left( {{\rm
{\bf a}}_{.s}^{(k)} }  \right)} \right) {\kern 1pt} _{\beta}
^{\beta} } } }}{{{\sum\limits_{\beta \in J_{r,\,\,n}} {{\left|
{\left( {{\rm {\bf A}}^{ k+1} } \right){\kern 1pt}
_{\beta} ^{\beta} }  \right|}}} }}}\cdot b_{sj}=
\]
\[
{\frac{{{\sum\limits_{\beta \in J_{r,\,n} {\left\{ {i}
\right\}}}\sum_{s=1}^{n} {{\rm{cdet}} _{i} \left( {\left( {{\rm
{\bf A}}^{ k+1} } \right)_{\,. \,i} \left( {{\rm {\bf
a}}_{.s}^{ (k)} }  \right)} \right) {\kern 1pt} _{\beta} ^{\beta} }
} }\cdot b_{sj}} {{{\sum\limits_{\beta \in J_{r,\,\,n}} {{\left|
{\left( {{\rm {\bf A}}^{ k+1} } \right){\kern 1pt}
_{\beta} ^{\beta} }  \right|}}} }}}
\]
Since ${\sum\limits_{s} {{\rm {\bf a}}_{.\,s}^{ (k)}  b_{sj}} }=
\left( {{\begin{array}{*{20}c}
 {{\sum\limits_{s} {a_{1s}^{(k)}  b_{sj}} } } \hfill \\
 {{\sum\limits_{k} {a_{2s}^{(k)}  b_{sj}} } } \hfill \\
 { \vdots}  \hfill \\
 {{\sum\limits_{k} {a_{ns}^{ (k)}  b_{sj}} } } \hfill \\
\end{array}} } \right) = \hat{{\rm {\bf b}}}_{.j}$, where $\hat{{\rm {\bf b}}}_{.j} $ denotes  the $j$th
column of $\hat{{\rm {\bf B}}}$ for all  $j = \overline {1,m} $, then it follows (\ref{eq:dr_AX}).
$\blacksquare$

For complex matrix equation (\ref{eq:AX=B}) we evident have the following corollaries, where ${\bf A}$ is not necessarily Hermitian.
 \begin{Corollary}(\cite{ky2}, Theorem 3.2.)  If $Ind{\kern 1pt} {\rm {\bf A}} = k $ and $\rank{\rm {\bf
A}}^{k + 1} = \rank{\rm {\bf A}}^{k}=r \le n$ for ${\rm {\bf A}}\in {\mathbb C}^{n\times n} $, then for Drazin inverse solution $ {\bf X}={\bf A}^{D}{\bf B}= (x_{ij})$ of  (\ref{eq:AX=B}) we have
\begin{equation*}
 x_{ij} = {\frac{{{\sum\limits_{\beta \in
J_{r,\,n} {\left\{ {i} \right\}}} {{\left| \left( {\left( {{\rm
{\bf A}}^{ k+1}} \right)_{\,.\,i} \left( {\hat{{\rm
{\bf b}}}_{.j}} \right)} \right){\kern 1pt} {\kern 1pt} _{\beta}
^{\beta} \right|}} } }}{{{\sum\limits_{\beta \in J_{r,\,\,n}}
{{\left| {\left( {{\rm {\bf A}}^{ k+1} } \right){\kern
1pt} {\kern 1pt} _{\beta} ^{\beta} }  \right|}}} }}}.
\end{equation*}
 \end{Corollary}
 \begin{Corollary}(\cite{ky1}, Theorem 4.5.)  If $Ind{\kern 1pt} {\rm {\bf A}} = k $ and $\rank{\rm {\bf
A}}^{k + 1} = \rank{\rm {\bf A}}^{k}=r \le n$ for ${\rm {\bf A}}\in {\mathbb C}^{n\times n} $, and ${\rm {\bf y}} = \left( {y_{1},\ldots ,y_{n} } \right)^{T}\in {\mathbb C}^{n}$, then for Drazin inverse solution $ {\bf x}={\bf A}^{D}{\bf y}=:  \left( {x_{1},\ldots ,x_{n} } \right)^{T}\in {\mathbb C}^{n}$ of the system of  linear equations
\[{\rm {\bf A}} \cdot {\rm {\bf x}} = {\rm {\bf y}},\]
we have for all $j = \overline {1,n}$,

\[
x_{j} = {\frac{{\sum\limits_{\beta \in J_{r,n}
{\left\{ {j} \right\}}} {{\left| {\left( {({\rm {\bf A}}^{k+1})_{.\,j} ({\rm {\bf f}}  )} \right)_{\beta} ^{\beta} }
\right|}}}}{{\left|  \left( {{\rm {\bf A}}^{k+1} }
\right)_{\beta} ^{\beta}
\right| } }},
\]
where  ${\rm {\bf f}} = {\rm {\bf A}}^{ k} {\rm {\bf y}}$.
\end{Corollary}

Consider a  matrix equation
\begin{equation}\label{eq:XA=B}
 {\rm {\bf X}}{\rm {\bf A}} = {\rm {\bf B}},
\end{equation}
 where $ {\bf
A}\in {\mathbb{H}}^{n\times n} $, $ {\rm {\bf B}}\in {\mathbb{H}}^{m\times n} $ are given, ${\rm {\bf
A}}$ is Hermitian and ${\rm {\bf X}} \in
{\mathbb{H}}^{m\times n}$ is unknown.
 Let $Ind{\kern 1pt} {\rm {\bf A}} = k $
and denote  ${\rm {\bf B}}{\rm {\bf A}}^{
k}=:\check{{\rm {\bf B}}}= (\check{b}_{ij})\in
{\mathbb{H}}^{m\times n}$.

\begin{theorem}
If $\rank{\rm {\bf
A}}^{k + 1} = \rank{\rm {\bf A}}^{k}=r \le n$ for ${\rm {\bf A}}\in {\mathbb H}^{n\times n} $, then for the Drazin inverse solution $ {\bf X}={\bf B}{\bf A}^{D}=: (x_{ij})$ of  (\ref{eq:XA=B}),  we have for all $i =
\overline {1,m}$, $j =
\overline {1,n} $
\begin{equation}
\label{eq:dr_XA} x_{ij} = {\frac{{{\sum\limits_{\alpha \in I_{r,n}
{\left\{ {j} \right\}}} {{\rm{rdet}} _{j} \left( {\left( {{\rm {\bf A}}^{k+1} } \right)_{\,j\,.} \left( {\check{{\rm
{\bf b}}}_{i\,.}} \right)} \right)\,_{\alpha} ^{\alpha} }
}}}{{{\sum\limits_{\alpha \in I_{r,\,n}}  {{\left| {\left( {{\rm {\bf A}}^{ k+1} } \right) {\kern 1pt} _{\alpha}
^{\alpha} } \right|}}} }}}.
\end{equation}

\noindent where $\check{{\rm {\bf b}}}_{i.}$ is the $i$th row of
$\check{{\rm {\bf B}}}$ for all $i =
\overline {1,m}$.

\end{theorem}
\textit{Proof}  By Theorem  \ref{theor:det_rep_draz} we can represent the matrix
${\rm {\bf A}}^{ D} $ by (\ref{eq:dr_rep_rdet}). Therefore, for
all for all $i =
\overline {1,m}$, $j =
\overline {1,n} $,  we obtain
 \begin{gather*}
x_{ij}  =\sum_{s=1}^{n}b_{is} a_{sj}^{D}=\sum_{s=1}^{n}b_{is}\cdot
{\frac{{{\sum\limits_{\alpha \in I_{r,n} {\left\{ {j} \right\}}}
{{\rm{rdet}} _{j} \left( {\left( {{\rm {\bf A}}^{k+1}
} \right)_{\,j\,.} \left( {\bf a}_{s\,.}^{(k)} \right)}
\right)\,_{\alpha} ^{\alpha} } }}}{{{\sum\limits_{\alpha \in
I_{r,\,n}}  {{\left| {\left( {{\rm {\bf A}}^{ k+1} }
\right) {\kern 1pt} _{\alpha} ^{\alpha} } \right|}}} }}} =
\\
{\frac{{{\sum_{s=1}^{n}b_{is} \sum\limits_{\alpha \in I_{r,n}
{\left\{ {j} \right\}}} {{\rm{rdet}} _{j} \left( {\left( {{\rm {\bf A}}^{ k+1} } \right)_{\,j\,.} \left( {\bf
a}_{s\,.}^{(k)} \right)} \right)\,_{\alpha} ^{\alpha} } } }}
{{{\sum\limits_{\alpha \in I_{r,\,n}}  {{\left| {\left( {{\rm {\bf A}}^{ k+1} } \right) {\kern 1pt} _{\alpha} ^{\alpha} }
\right|}}} }}}
 \end{gather*}
Since ${\sum\limits_{s} {{  b_{is}\rm {\bf a}}_{s\,.}^{(k)}}
}=\begin{pmatrix}
  \sum\limits_{s} {b_{is}a_{s1}^{(k)} } & \sum\limits_{s} {b_{is}a_{s2}^{(k)} } & \cdots & \sum\limits_{s} {b_{is}a_{sn}^{(k)} }
\end{pmatrix}
 = \check{{\rm {\bf b}}}_{i.}$, where $\check{{\rm {\bf b}}}_{i.}$ denotes the $i$th row of
$\check{{\rm {\bf B}}}$ for all $i =
\overline {1,m}$,
then it follows (\ref{eq:dr_XA}).
$\blacksquare$

We evident have the following corollary for the complex matrix equation (\ref{eq:XA=B}), where ${\bf A}$ is not necessarily Hermitian.
\begin{Corollary}(\cite{ky2}, Theorem 3.4.) If $\rank{\rm {\bf
A}}^{k + 1} = \rank{\rm {\bf A}}^{k}=r \le n$ for ${\rm {\bf A}}\in {\mathbb C}^{n\times n} $, then for the Drazin inverse solution $ {\bf X}={\bf B}{\bf A}^{D}=: (x_{ij})$ of  (\ref{eq:XA=B}),  we have for all $i =
\overline {1,m}$, $j =
\overline {1,n}, $
\begin{equation*}
 x_{ij} = {\frac{{{\sum\limits_{\alpha \in I_{r,n}
{\left\{ {j} \right\}}} {{\left| \left( {\left( {{\rm
{\bf A}}^{ k+1} } \right)_{\,j\,.} \left( {\check{{\rm {\bf
b}}}_{i\,.}} \right)} \right)\,_{\alpha} ^{\alpha}\right|} }
}}}{{{\sum\limits_{\alpha \in I_{r,\,n}}  {{\left| {\left( {{\rm {\bf A}}^{ k+1} } \right) {\kern 1pt} _{\alpha}
^{\alpha} } \right|}}} }}}.
\end{equation*}
\end{Corollary}

Consider a  matrix equation
 \begin{equation}\label{eq1:AXB=D}
 {\rm {\bf A}}{\rm {\bf X}}{\rm {\bf B}} = {\rm {\bf
D}},
\end{equation}
where $ {\rm {\bf
A}}\in{\rm {\mathbb{H}}}^{n \times n}$, $ {\rm {\bf B}}\in{\rm {\mathbb{H}}}^{m \times m}$, $ {\rm {\bf D}}\in{\rm {\mathbb{H}}}^{n \times m}$ are given, ${\rm {\bf
A}}$, ${\rm {\bf B}}$ are Hermitian, and $ {\rm {\bf X}}\in{\rm
{\mathbb{H}}}^{n \times m}$ is unknown.
Let $Ind{\kern 1pt} {\rm {\bf A}} = k_{1} $ and $Ind{\kern 1pt} {\rm {\bf B}} = k_{2} $ and denote  ${\rm {\bf A}}^{
 k_{1} }{\rm {\bf D}}{\rm {\bf B}}^{
 k_{2} }=:\widetilde{{\rm {\bf D}}}= (\widetilde{d}_{ij})\in
{\mathbb{H}}^{n\times m}$.

\begin{theorem}\label{theor:AXB=D}
If  $\rank{\rm {\bf
A}}^{k_{1} + 1} = \rank{\rm {\bf A}}^{k_{1}}=r_{1} \le n$ for ${\rm {\bf A}}\in {\mathbb H}^{n\times n} $, and  $\rank{\rm {\bf
B}}^{k_{2} + 1} = \rank{\rm {\bf B}}^{k_{2}}=r_{2} \le m$ for ${\rm {\bf B}}\in {\mathbb H}^{m\times m} $, then for the
Drazin inverse solution ${\rm {\bf X}}={\rm {\bf A}}^{D}{\rm {\bf D}}{\rm {\bf B}}^{
D}=(x_{ij})\in
{\mathbb{H}}^{n\times m}$  of (\ref{eq1:AXB=D}) we have
\begin{equation}\label{eq:d^B}
x_{ij} = {\frac{{{\sum\limits_{\beta \in J_{r_{1},\,n} {\left\{
{i} \right\}}} {{\rm{cdet}} _{i} \left( {\left( {{\rm {\bf A}}^{{k_{1}}+1
}} \right)_{\,.\,i} \left( {{{\rm {\bf
d}}}\,_{.\,j}^{{\rm {\bf B}}}} \right)} \right) _{\beta} ^{\beta}
} } }}{{{\sum\limits_{\beta \in J_{r_{1},n}} {{\left| {\left(
{{\rm {\bf A}}^{{k_{1}}+1}} \right)_{\beta} ^{\beta} }
\right|}} \sum\limits_{\alpha \in I_{r_{2},m}}{{\left| {\left(
{{\rm {\bf B}}^{{k_{2}}+1
} } \right) _{\alpha} ^{\alpha} }
\right|}}} }}},
\end{equation}
or
\begin{equation}\label{eq:d^A}
 x_{ij}={\frac{{{\sum\limits_{\alpha
\in I_{r_{2},m} {\left\{ {j} \right\}}} {{\rm{rdet}} _{j} \left(
{\left( {{\rm {\bf B}}^{{k_{2}}+1
} } \right)_{\,j\,.} \left(
{{{\rm {\bf d}}}\,_{i\,.}^{{\rm {\bf A}}}} \right)}
\right)\,_{\alpha} ^{\alpha} } }}}{{{\sum\limits_{\beta \in
J_{r_{1},n}} {{\left| {\left( {{\rm {\bf A}}^{{k_{1}}+1
}}
\right) _{\beta} ^{\beta} } \right|}}\sum\limits_{\alpha \in
I_{r_{2},m}} {{\left| {\left( {{\rm {\bf B}}^{{k_{2}}+1
} }
\right) _{\alpha} ^{\alpha} } \right|}}} }}},
\end{equation}
where
\begin{equation} \label{eq:def_d^B_m}
   {{{\rm {\bf d}}}_{.\,j}^{{\rm {\bf B}}}}=\left(
\sum\limits_{\alpha \in I_{r_{2},m} {\left\{ {j} \right\}}}
{{\rm{rdet}} _{j} \left( {\left( {{\rm {\bf B}}^{{k_{1}}+1}
} \right)_{j.} \left( {\tilde{{\rm {\bf d}}}_{l.}} \right)}
\right)_{\alpha} ^{\alpha}}
\right)\in{\rm {\mathbb{H}}}^{n \times
1},\,\,\,\,l=1,...,n \end{equation}
\begin{equation} \label{eq:d_A<n}  {{{\rm {\bf d}}}_{i\,.}^{{\rm {\bf A}}}}=\left(
\sum\limits_{\beta \in J_{r_{1},n} {\left\{ {i} \right\}}}
{{\rm{cdet}} _{i} \left( {\left( {{\rm {\bf A}}^{{k_{1}}+1}
} \right)_{.i} \left( {\tilde{{\rm {\bf d}}}_{.t}} \right)}
\right)_{\beta} ^{\beta}} \right)\in{\rm {\mathbb{H}}}^{1 \times
m},\,\,\,\,t=1,...,m
\end{equation}
 are the column vector and the row vector, respectively.  ${\tilde{{\rm {\bf
 d}}}_{i.}}$ and
${\tilde{{\rm {\bf d}}}_{.j}}$ are the ith row   and the jth
column  of ${\rm {\bf \widetilde{D}}}$ for all $i =
\overline {1,n}$, $j =
\overline {1,m}$.

\end{theorem}
\textit{Proof.}
An entry of the Drazin inverse solution ${\rm {\bf X}}={\rm {\bf A}}^{D}{\rm {\bf D}}{\rm {\bf B}}^{
D}=(x_{ij})\in
{\mathbb{H}}^{n\times m}$ is

\begin{equation}
\label{eq:sum+} x_{ij} = {{\sum\limits_{s = 1}^{m} {\left(
{{\sum\limits_{t = 1}^{n} {{a}_{it}^{D} d_{ts}} } } \right)}}
{b}_{sj}^{D}}
\end{equation}
 for all $i =
\overline {1,n}$, $j =
\overline {1,m}$, where by Theorem \ref{theor:det_rep_draz}

\[
 a_{ij}^{ D}  = {\frac{{{\sum\limits_{\beta
\in J_{r_{1},\,n} {\left\{ {i} \right\}}} {{\rm{cdet}} _{i} \left(
{\left( {{\rm {\bf A}}^{k_{1}+1} } \right)_{\,. \,i}
\left( {{\rm {\bf a}}_{.j}^{( k_{1})} }  \right)} \right){\kern 1pt}
{\kern 1pt} _{\beta} ^{\beta} } } }}{{{\sum\limits_{\beta \in
J_{r_{1},\,n}} {{\left| {\left( {{\rm {\bf A}}^{ {k_{1}+1}}}
\right){\kern 1pt} _{\beta} ^{\beta} }  \right|}}} }}},
\]
\begin{equation}\label{eq:b+}
 b_{ij}^{ D}  =
{\frac{{{\sum\limits_{\alpha \in I_{r_{2},m} {\left\{ {j}
\right\}}} {{\rm{rdet}} _{j} \left( {({\rm {\bf B}}^{
k_{2}+1} )_{j\,.\,} ({\rm {\bf b}}_{i.\,}^{ (k_{2})} )} \right)\,_{\alpha}
^{\alpha} } }}}{{{\sum\limits_{\alpha \in I_{r_{2},m}}  {{\left|
{\left( {{\rm {\bf B}}^{
k_{2}+1} } \right){\kern 1pt}
_{\alpha} ^{\alpha} } \right|}}} }}}.
\end{equation}

Denote  by $\hat{{\rm {\bf d}}_{.s}}$ the $s$th column of ${\rm
{\bf A}}^{ k_{1}}{\rm {\bf D}}=:\hat{{\rm {\bf D}}}=
(\hat{d}_{ij})\in {\mathbb{H}}^{n\times m}$ for all $s=\overline
{1,m}$. It follows from ${\sum\limits_{l} { {\rm {\bf a}}_{.\,l}^{
 (k_{1})}}d_{ls} }=\hat{{\rm {\bf d}}_{.\,s}}$ that
\[
\sum\limits_{l = 1}^{n} {{a}_{il}^{D} d_{ls}}=\sum\limits_{l =
1}^{n}{\frac{{{\sum\limits_{\beta \in J_{r_{1},\,n} {\left\{ {i}
\right\}}} {{\rm{cdet}} _{i} \left( {\left( {{\rm {\bf A}}^{ k_{1}+1}
} \right)_{\,. \,i} \left( {{\rm {\bf a}}_{.l}^{(k_{1})} }
\right)} \right){\kern 1pt} {\kern 1pt} _{\beta} ^{\beta} } }
}}{{{\sum\limits_{\beta \in J_{r_{1},\,n}} {{\left| {\left( {{\rm
{\bf A}}^{ k_{1}+1}} \right){\kern 1pt} _{\beta} ^{\beta}
}  \right|}}} }}}\cdot d_{ls}=
\]
\begin{equation}\label{eq:sum_cdet}
{\frac{{{\sum\limits_{\beta \in J_{r_{1},\,n} {\left\{ {i}
\right\}}}\sum\limits_{l = 1}^{n} {{\rm{cdet}} _{i} \left( {\left(
{{\rm {\bf A}}^{ k_{1}+1}} \right)_{\,. \,i} \left( {{\rm
{\bf a}}_{.l}^{ (k_{1})} } \right)} \right){\kern 1pt} {\kern 1pt}
_{\beta} ^{\beta} } } }\cdot d_{ls}}{{{\sum\limits_{\beta \in
J_{r_{1},\,n}} {{\left| {\left( {{\rm {\bf A}}^{ k_{1}+1}}
\right){\kern 1pt} _{\beta} ^{\beta} }  \right|}}}
}}}={\frac{{{\sum\limits_{\beta \in J_{r_{1},\,n} {\left\{ {i}
\right\}}} {{\rm{cdet}} _{i} \left( {\left( {{\rm {\bf A}}^{ k_{1}+1}
} \right)_{\,. \,i} \left( \hat{{\rm {\bf d}}_{.\,s}}
\right)} \right){\kern 1pt} {\kern 1pt} _{\beta} ^{\beta} } }
}}{{{\sum\limits_{\beta \in J_{r_{1},\,n}} {{\left| {\left( {{\rm
{\bf A}}^{ k_{1}+1}} \right){\kern 1pt} _{\beta} ^{\beta}
}  \right|}}} }}}
\end{equation}
Suppose ${\rm {\bf e}}_{s.}$ and ${\rm {\bf e}}_{.\,s}$ are
respectively the unit row-vector and the unit column-vector whose
components are $0$, except the $s$th components, which are $1$.
Substituting  (\ref{eq:sum_cdet}) and (\ref{eq:b+}) in
(\ref{eq:sum+}), we obtain
\[
x_{ij} =\sum\limits_{s = 1}^{m}{\frac{{{\sum\limits_{\beta \in
J_{r_{1},\,n} {\left\{ {i} \right\}}} {{\rm{cdet}} _{i} \left(
{\left( {{\rm {\bf A}}^{ k_{1}+1} } \right)_{\,. \,i}
\left( \hat{{\rm {\bf d}}_{.\,s}} \right)} \right){\kern 1pt}
{\kern 1pt} _{\beta} ^{\beta} } } }}{{{\sum\limits_{\beta \in
J_{r_{1},\,n}} {{\left| {\left( {{\rm {\bf A}}^{ k_{1}+1}}
\right){\kern 1pt} _{\beta} ^{\beta} }  \right|}}}
}}}{\frac{{{\sum\limits_{\alpha \in I_{r_{2},m} {\left\{ {j}
\right\}}} {{\rm{rdet}} _{j} \left( {({\rm {\bf B}}^{ k_{2}+1} )_{j\,.\,} ({\rm {\bf b}}_{s.\,}^{ (k_{2})} )} \right)\,_{\alpha}
^{\alpha} } }}}{{{\sum\limits_{\alpha \in I_{r_{2},m}}  {{\left|
{\left( {{\rm {\bf B}}^{ k_{2}+1} } \right){\kern 1pt}
_{\alpha} ^{\alpha} } \right|}}} }}}.
\]
Since \begin{equation}\label{eq:prop}\hat{{\rm {\bf
d}}_{.\,s}}=\sum\limits_{l = 1}^{n}{\rm {\bf e}}_{.\,l}\hat{
d_{ls}},\,  {\rm {\bf b}}_{s.\,}^{(k_{2})}=\sum\limits_{t =
1}^{m}b_{st}^{(k_{2})}{\rm {\bf
e}}_{t.},\,\sum\limits_{s=1}^{m}\hat{d_{ls}}b_{st}^{(k_{2})}=\widetilde{d}_{lt},\end{equation}
then we have
\[
x_{ij} = \]
\[{\frac{{ \sum\limits_{s = 1}^{m}\sum\limits_{t =
1}^{m} \sum\limits_{l = 1}^{n} {\sum\limits_{\beta \in
J_{r_{1},\,n} {\left\{ {i} \right\}}} {{\rm{cdet}} _{i} \left(
{\left( {{\rm {\bf A}}^{ k_{1}+1}} \right)_{\,. \,i}
\left( {\rm {\bf e}}_{.\,l} \right)} \right){\kern 1pt} {\kern
1pt} _{\beta} ^{\beta} } } }\hat{
d_{ls}}b_{st}^{ (k_{2})}{\sum\limits_{\alpha \in I_{r_{2},m} {\left\{ {j}
\right\}}} {{\rm{rdet}} _{j} \left( {({\rm {\bf B}}^{ k_{2}+1})_{j\,.\,} ({\rm {\bf e}}_{t.} )} \right)\,_{\alpha} ^{\alpha}
} } }{{{\sum\limits_{\beta \in J_{r_{1},\,n}} {{\left| {\left(
{{\rm {\bf A}}^{ k_{1}+1}} \right){\kern 1pt} _{\beta}
^{\beta} }  \right|}}} }{{\sum\limits_{\alpha \in I_{r_{2},m}}
{{\left| {\left( {{\rm {\bf B}}^{ k_{2}+1} } \right){\kern
1pt} _{\alpha} ^{\alpha} } \right|}}} }}    }=
\]
\begin{equation}\label{eq:x_ij}
{\frac{{ \sum\limits_{t = 1}^{m} \sum\limits_{l = 1}^{n}
{\sum\limits_{\beta \in J_{r_{1},\,n} {\left\{ {i} \right\}}}
{{\rm{cdet}} _{i} \left( {\left( {{\rm {\bf A}}^{ k_{1}+1}} \right)_{\,. \,i} \left( {\rm {\bf e}}_{.\,l} \right)}
\right){\kern 1pt} {\kern 1pt} _{\beta} ^{\beta} } }
}\,\,\widetilde{d}_{lt}{\sum\limits_{\alpha \in I_{r_{2},m}
{\left\{ {j} \right\}}} {{\rm{rdet}} _{j} \left( {({\rm {\bf B}}^{ k_{2}+1} )_{j\,.\,} ({\rm {\bf e}}_{t.} )}
\right)\,_{\alpha} ^{\alpha} } } }{{{\sum\limits_{\beta \in
J_{r_{1},\,n}} {{\left| {\left( {{\rm {\bf A}}^{ k_{1}+1}}
\right){\kern 1pt} _{\beta} ^{\beta} }  \right|}}}
}{{\sum\limits_{\alpha \in I_{r_{2},m}} {{\left| {\left( {{\rm {\bf B}}^{ k_{2}+1}} \right){\kern 1pt} _{\alpha}
^{\alpha} } \right|}}} }}    }.
\end{equation}
Denote by
\[
 d^{{\rm {\bf A}}}_{it}:= \]
\[
{\sum\limits_{\beta \in J_{r_{1},\,n} {\left\{ {i} \right\}}}
{{\rm{cdet}} _{i} \left( {\left( {{\rm {\bf A}}^{ k_{1}+1}} \right)_{\,. \,i} \left( \widetilde{{\rm {\bf d}}}_{.\,t}
\right)} \right) _{\beta} ^{\beta} } }= \sum\limits_{l
= 1}^{n} {\sum\limits_{\beta \in J_{r_{1},\,n} {\left\{ {i}
\right\}}} {{\rm{cdet}} _{i} \left( {\left( {{\rm {\bf A}}^{ k_{1}+1}} \right)_{\,. \,i} \left( {\rm {\bf e}}_{.\,l}
\right)} \right)_{\beta} ^{\beta} } }{\kern 1pt}
\widetilde{d}_{lt}
\]
the $t$th component  of a row-vector ${\rm {\bf d}}^{{\rm {\bf
A}}}_{i\,.}= (d^{{\rm {\bf A}}}_{i1},...,d^{{\rm {\bf A}}}_{im})$
for all  $t=
\overline {1,m}$. Substituting it in (\ref{eq:x_ij}),
we have
\[x_{ij} ={\frac{{ \sum\limits_{t = 1}^{m}
 d^{{\rm {\bf A}}}_{it}
}{\sum\limits_{\alpha \in I_{r_{2},m} {\left\{ {j} \right\}}}
{{\rm{rdet}} _{j} \left( {({\rm {\bf B}}^{ k_{2}+1}
)_{j\,.\,} ({\rm {\bf e}}_{t.} )} \right)_{\alpha} ^{\alpha} } }
}{{{\sum\limits_{\beta \in J_{r_{1},\,n}} {{\left| {\left( {{\rm
{\bf A}}^{ k_{1}+1}} \right){\kern 1pt} _{\beta} ^{\beta}
}  \right|}}} }{{\sum\limits_{\alpha \in I_{r_{2},m}} {{\left|
{\left( {{\rm {\bf B}}^{ k_{2}+1} } \right)
_{\alpha} ^{\alpha} } \right|}}} }}    }.
\]
Since $\sum\limits_{t = 1}^{m}
 d^{{\rm {\bf A}}}_{it}{\rm {\bf e}}_{t.}={\rm {\bf
d}}^{{\rm {\bf A}}}_{i\,.}$, then it follows (\ref{eq:d^A}).

If we denote by
\begin{equation}\label{eq:d^B_den}
 d^{{\rm {\bf B}}}_{lj}:=
\sum\limits_{t = 1}^{m}\widetilde{d}_{lt}{\sum\limits_{\alpha \in
I_{r_{2},m} {\left\{ {j} \right\}}} {{\rm{rdet}} _{j} \left(
{({\rm {\bf B}}^{ k_{2}+1} )_{j\,.} ({\rm {\bf e}}_{t.}
)} \right)_{\alpha} ^{\alpha} } }={\sum\limits_{\alpha \in
I_{r_{2},m} {\left\{ {j} \right\}}} {{\rm{rdet}} _{j} \left(
{({\rm {\bf B}}^{ k_{2}+1} )_{j\,.\,} (\widetilde{{\rm {\bf
d}}}_{l.} )} \right)_{\alpha} ^{\alpha} } }
\end{equation}

\noindent the $l$th component  of a column-vector ${\rm {\bf
d}}^{{\rm {\bf B}}}_{.\,j}= (d^{{\rm {\bf B}}}_{1j},...,d^{{\rm
{\bf B}}}_{nj})^{T}$ for all $l=1,...,n$ and substitute it
in (\ref{eq:x_ij}), we obtain
\[x_{ij} ={\frac{{  \sum\limits_{l = 1}^{n}
{\sum\limits_{\beta \in J_{r_{1},n} {\left\{ {i} \right\}}}
{{\rm{cdet}} _{i} \left( {\left( {{\rm {\bf A}}^{ k_{1}+1}} \right)_{\,. \,i} \left( {\rm {\bf e}}_{.\,l} \right)}
\right) _{\beta} ^{\beta} } } }\,\,d^{{\rm {\bf
B}}}_{lj} }{{{\sum\limits_{\beta \in J_{r_{1},\,n}} {{\left|
{\left( {{\rm {\bf A}}^{ k_{1}+1}} \right)
_{\beta} ^{\beta} }  \right|}}} }{{\sum\limits_{\alpha \in
I_{r_{2},m}} {{\left| {\left( {{\rm {\bf B}}^{ k_{2}+1} }
\right) _{\alpha} ^{\alpha} } \right|}}} }}    }.
\]
Since $\sum\limits_{l = 1}^{n}{\rm {\bf e}}_{.l}
 d^{{\rm {\bf B}}}_{lj}={\rm {\bf
d}}^{{\rm {\bf B}}}_{.j}$, then it follows (\ref{eq:d^B}).
$\blacksquare$

\begin{Corollary}(\cite{ky2}, Theorem 3.6.)
 If  $\rank{\rm {\bf
A}}^{k_{1} + 1} = \rank{\rm {\bf A}}^{k_{1}}=r_{1} \le n$ for ${\rm {\bf A}}\in {\mathbb C}^{n\times n} $, and  $\rank{\rm {\bf
B}}^{k_{2} + 1} = \rank{\rm {\bf B}}^{k_{2}}=r_{2} \le m$ for ${\rm {\bf B}}\in {\mathbb C}^{m\times m} $, then for the
Drazin inverse solution ${\rm {\bf X}}={\rm {\bf A}}^{D}{\rm {\bf D}}{\rm {\bf B}}^{
D}=:(x_{ij})\in
{\mathbb{C}}^{n\times m}$  of (\ref{eq1:AXB=D}) we have
\[
x_{ij} = {\frac{{{\sum\limits_{\beta \in J_{r_{1},\,n} {\left\{
{i} \right\}}} { \left| {{\rm {\bf A}}^{ k_{1}+1} _{\,.\,i} \left( {{{\rm {\bf d}}}\,_{.\,j}^{{\rm {\bf B}}}}
\right)\, _{\beta} ^{\beta}} \right| } } }}{{{\sum\limits_{\beta
\in J_{r_{1},n}} {{\left| {\left( {{\rm {\bf A}}^{ k_{1}+1} } \right)_{\beta} ^{\beta} } \right|}} \sum\limits_{\alpha \in
I_{r_{2},m}}{{\left| {\left( {{\rm {\bf B}}^{k_{2}+1} }
\right) _{\alpha} ^{\alpha} } \right|}}} }}},
\]
or
\[
 x_{ij}={\frac{{{\sum\limits_{\alpha
\in I_{r_{2},m} {\left\{ {j} \right\}}} { \left| {{\rm {\bf B}}^{k_{2}+1} _{j\,.} \left( {{{\rm {\bf
d}}}\,_{i\,.}^{{\rm {\bf A}}}} \right)\,_{\alpha} ^{\alpha}}
\right| } }}}{{{\sum\limits_{\beta \in J_{r_{1},n}} {{\left|
{\left( {{\rm {\bf A}}^{k_{1}+1} } \right) _{\beta}
^{\beta} } \right|}}\sum\limits_{\alpha \in I_{r_{2},m}} {{\left|
{\left( {{\rm {\bf B}}^{k_{2}+1} } \right) _{\alpha}
^{\alpha} } \right|}}} }}},
\]
where
  \[
   {{\rm {\bf d}}_{.\,j}^{{\rm {\bf B}}}}=\left[
\sum\limits_{\alpha \in I_{r_{2},m} {\left\{ {j} \right\}}} {
\left| {{\rm {\bf B}}^{ k_{2}+1} _{j.}
\left( {\widetilde{{\rm {\bf d}}}_{1.}} \right)\,_{\alpha} ^{\alpha}}
\right|},...,\sum\limits_{\alpha \in I_{r_{2},m} {\left\{ {j}
\right\}}} { \left| {{\rm {\bf B}}^{k_{2}+1}_{j.} \left( {\widetilde{{\rm {\bf d}}}_{n.}}
\right)\,_{\alpha} ^{\alpha}} \right|} \right]^{T},
\]
  \[
  {{\rm {\bf d}}_{i\,.}^{{\rm {\bf A}}}}=\left[
\sum\limits_{\beta \in J_{r_{1},n} {\left\{ {i} \right\}}} {
\left| {{\rm {\bf A}}^{k_{1}+1}_{.i}
\left( {\widetilde{{\rm {\bf d}}}_{.1}} \right)\,_{\beta} ^{\beta}}
\right|},...,\sum\limits_{\alpha \in I_{r_{1},n} {\left\{ {i}
\right\}}} { \left| {{\rm {\bf A}}^{k_{1}+1}_{.i} \left( {\widetilde{{\rm {\bf d}}}_{.\,m}}
\right)\,_{\beta} ^{\beta}} \right|} \right]
\]
 are the column-vector and the row-vector. ${\widetilde{{\rm {\bf d}}}_{i.}}$  and ${\widetilde{{\rm {\bf d}}}_{.j}}$ are respectively the $i$th row  and the $j$th column of $\widetilde{{\rm {\bf D}}}$ for all $i=\overline {1,n}$, $j=\overline {1,m}$.

\end{Corollary}

\section{An example}
In this section, we give an example to illustrate our results. Let
us consider the matrix equation
\begin{equation}\label{eq_ex:AXB=D}
 {\rm {\bf A}}{\rm {\bf X}}{\rm {\bf B}} = {\rm {\bf
D}},
\end{equation}
where
\[{\bf A}=\begin{pmatrix}
  1 & k & -i \\
  -k & 2 & j \\
 i & -j & 1
\end{pmatrix},\,\, {\bf B}=\begin{pmatrix}
 1 & i \\
 -i & 1
\end{pmatrix},\,\, {\bf D}=\begin{pmatrix}
  1 & i \\
  k & 1 \\
  1 & j
\end{pmatrix}.\]
Since ${\bf A}^{2}=\begin{pmatrix}
  3 & 4k & -3i \\
  -4k & 6 & 4j \\
 3i & -4j & 3
\end{pmatrix}$,
 $\det {\bf A}=\det {\bf A}^{2}=0$, and $\det \begin{pmatrix}
  1 & k  \\
  -k & 2
\end{pmatrix}=1$, $\det \begin{pmatrix}
  3 & 4k \\
  -4k & 6
\end{pmatrix}=2$, then, by Theorem \ref{theor:rank_her}, $Ind\,{\bf A}=1$ and $r_{1}=\rank {\bf A}=2$.  Similarly, since $ {\bf B}^{2}=\begin{pmatrix}
 2 &2 i \\
 -2i & 2
\end{pmatrix}$, then $Ind\,{\bf B}=1$ and $r_{2}=\rank {\bf B}=1$.

 We shall
find the Drazin inverse solution ${\rm {\bf X}}^{d}=(x_{ij}^{d})$ of (\ref{eq_ex:AXB=D}) by the equations
(\ref{eq:d^B})-(\ref{eq:def_d^B_m}). We obtain \[ {{{\sum\limits_{\alpha \in I_{1,\,2}}
{{\left| {\left( {{\rm {\bf B}}^{2} } \right) {\kern
1pt} _{\alpha} ^{\alpha} } \right|}}} }}=2+2=4,\]
\[ {{{\sum\limits_{\beta \in J_{2,\,3}} {{\left| {\left( {{\rm {\bf
A}}^{2}} \right) {\kern 1pt} _{\beta} ^{\beta} }
\right|}}} }}=
  \det\begin{pmatrix}
  3 & 4k \\
  -4k & 6
\end{pmatrix}+\det\begin{pmatrix}
  3 & -3i
 \\
3i & 3
\end{pmatrix}+\det\begin{pmatrix}
6 & 4j \\
 - 4j & 3
\end{pmatrix}=4.\]
Since
\[ {\rm {\bf \widetilde{D}}}= {\rm {\bf
A}}{\rm {\bf D}}{\rm {\bf B}}=\begin{pmatrix}
  1-i &1+i \\
 -i+j & 1-k \\
  1+i &  -1+i
\end{pmatrix}, \] then by (\ref{eq:def_d^B_m})
\[
   {{{\rm {\bf d}}}_{.\,j}^{{\rm {\bf B}}}}=\left(
\sum\limits_{\alpha \in I_{1,2} {\left\{ {j} \right\}}}
{{\rm{rdet}} _{j} \left( {\left( {{\rm {\bf B}}^{2}
} \right)_{1.} \left( {\tilde{{\rm {\bf d}}}_{l.}} \right)}
\right)_{\alpha} ^{\alpha}}
\right)\in{\rm {\mathbb{H}}}^{n \times
1},\,\,\,\,l=1,2,3\,\,\,j=1,2, \] and thus we have
\[{\bf d}_{.1}^{{\bf
B}}=\begin{pmatrix}
 1-i \\
-i+j \\
  1+i
\end{pmatrix},\,\,\,\,\,{\bf d}_{.2}^{{\bf  B}}=\begin{pmatrix}
1+i \\
 1-k \\
 -1+i
\end{pmatrix}.\]
Since
\[\left( {{\rm {\bf A}}^{ 2}} \right)_{\,.\,1}
\left( {{{\rm {\bf d}}}\,_{.\,1}^{{\rm {\bf B}}}} \right)=
\begin{pmatrix}
 1-i & 4k & -3i \\
  -i+j & 6 & 4j \\
  1+i & -4j & 3
\end{pmatrix},\] then finally we obtain
\[
  x_{11}^d = {\frac{{{\sum\limits_{\beta \in J_{2,\,3} {\left\{ {1}
\right\}}} { {\rm cdet} {\left(\left( {{\rm {\bf A}}^{ 2} } \right)_{\,.\,1} \left( {{{\rm {\bf d}}}\,_{.\,1}^{{\rm
{\bf B}}}} \right)\right)\, _{\beta} ^{\beta}}  } }
}}{{{\sum\limits_{\beta \in J_{2,3}} {{\left| {\left( {{\rm {\bf
A}}^{ 2} } \right)_{\beta} ^{\beta} } \right|}}
\sum\limits_{\alpha \in I_{1,2}}{{\left| {\left( {{\rm {\bf B}}^{ 2} } \right) _{\alpha} ^{\alpha} } \right|}}}
}}}=\]
 \[ \frac{{\rm cdet}_{1}\begin{pmatrix}
  1-i& 4k \\
 -i+j & 6
\end{pmatrix}+{\rm cdet}_{1}\begin{pmatrix}
  1-i & -3i \\
  1+i & 3
\end{pmatrix}}{16}=\frac{3-i+2j}{8}.
\]
Similarly,
\[
  x_{12}^d = \frac{{\rm cdet}_{1}\begin{pmatrix}
  1+i & 4k \\
  1-k & 6
\end{pmatrix}+{\rm cdet}_{1}\begin{pmatrix}
 1+i & -3i\\
  1+i  & 3
\end{pmatrix}}{16}=\frac{1+3i-2k}{8},\]
\[x_{21}^d =\frac{{\rm cdet}_{2}\begin{pmatrix}
 3 &  1-i \\
  -4k & -i+j
\end{pmatrix}+{\rm cdet}_{1}\begin{pmatrix}
  -i+j& 4j \\
 1+i &3
\end{pmatrix}}{16}=\frac{-3i-j+4k}{8},\]
 \[
   x_{22}^d =\frac{{\rm cdet}_{2}\begin{pmatrix}
 3 &  1+i \\
  -4k & 1-k
\end{pmatrix}+{\rm cdet}_{1}\begin{pmatrix}
  1-k& 4j \\
 -1+i &3
\end{pmatrix}}{16}=\frac{3+4j+k}{8},\]
\[   x_{31}^d =\frac{{\rm cdet}_{2}\begin{pmatrix}
  3 & 1-i \\
 3i& 1+i
\end{pmatrix}+{\rm cdet}_{2}\begin{pmatrix}
   6 & -i+j\\
  -4j & 1+i
\end{pmatrix}}{16}=\frac{1+3i+2k}{8},\]
\[
    x_{32}^d =\frac{{\rm cdet}_{2}\begin{pmatrix}
  3 & 1+i \\
 3i& -1+i
\end{pmatrix}+{\rm cdet}_{2}\begin{pmatrix}
   6 & 1-k\\
  -4j & -1+i
\end{pmatrix}}{16}=\frac{-3+i+2j}{8},\]
Thus,
\[{\rm {\bf X}}^d=\frac{1}{8}\left(
                                   \begin{array}{cc}
                                     3-i+2j & 1+3i-2k \\
                                    -3 i-j+4k & 3+4j+k \\
                                     1+3i+2k & -3+i+2j \\
                                   \end{array}
                                 \right)\]
is the Drazin inverse  solution of (\ref{eq_ex:AXB=D}).


\begin{thebibliography}{40}
\bibitem{dr}   M.P. Drazin, Pseudoinverse in associative rings and semigroups, Am. Math. Monthly 65 (1958) 506-514.
\bibitem{ca1} S. L. Campbell, C.D. Meyer, Generalized inverse of linear transformations, Corrected reprint of the 1979 original. Dover Publications, Inc., New York, 1991.
\bibitem{be} A. Ben-Israel, T. N. E. Greville, Generalized inverses: Theory and Applications, second ed., Springer, New York, 2003.
\bibitem{zh1}L. Zhang, A characterization of the Drazin inverse, Linear Algebra Appl. 335 (2001) 183-188.
\bibitem{ha} R.E.Hartwig, G. Wang, Y.M. Wei, Some additive results on Drazin inverse, Appl. Math. Comput. 322  (2001) 207-217.
\bibitem{kh} I.A. Khan, Q.W. Wang, The Drazin inverses in an arbitrary semiring, Linear and Multilinear Algebra 59 (9) (2011) 1019-1029.

 \bibitem{st} P.S. Stanimirovic', D.S. Djordjevic', Full-rank and determinantal representation
of the Drazin inverse, Linear Algebra Appl. 311 (2000) 131-
151.

    \bibitem{ky1} I.Kyrchei, Analogs of the adjoint matrix for generalized inverses and corresponding Cramer rules, Linear and Multilinear algebra 56(4) (2008)  453-469.
        \bibitem{ky11} I. Kyrchei, Analogs of Cramer's rule for the minimum norm least squares solutions of some matrix equations, Appl. Math. Comput. 218 (2012)  6375--6384.
 \bibitem{ky2}I. Kyrchei,
 Explicit formulas for determinantal representations of the Drazin inverse solutions of some matrix and differential matrix equations, Appl. Math. Comput. 219 (2013) 7632-7644.

 \bibitem{ky3} I. Kyrchei, Cramer's rule for quaternion
 systems of linear equations, Journal of Mathematical Sciences 155 (6) (2008) 839-858.

 \bibitem{ky4}   I. Kyrchei, The theory of the column and row determinants in a quaternion linear algebra,  in: Albert R. Baswell (Eds.), Advances in Mathematics Research 15,  Nova Sci. Publ., New York, 2012, pp. 301-359.

\bibitem{ky5} I. Kyrchei, Determinantal representations of the Moore-Penrose inverse over the quaternion skew field and corresponding Cramer's rules, Linear Multilinear Algebra, 59 (2011) 413-431.
\bibitem{ky33} I. Kyrchei, Determinantal representation of the Moore-Penrose inverse matrix over the quaternion skew field, Journal of Mathematical Sciences 180(1) (2012)  23--33.
\bibitem{ky6}    I. Kyrchei, Explicit representation formulas for the minimum norm least squares solutions of some quaternion matrix equations,  Linear Algebra Appl. 438 (2013) 136-152.
    \bibitem{so2} G. Song, Q.W. Wang, Condensed Cramer rule for some restricted
quaternion linear equations, Appl. Math. Comput. 218(7) (2011)
3110-3121.
 \bibitem{so1} G. Song, Q. Wang, H. Chang, Cramer rule for the unique solution of restricted matrix equations over the quaternion skew field, Comput Math. Appl. 61 (2011) 1576-1589.
\bibitem{song} G. Song, Determinantal representation of the generalized inverses over
the quaternion skew field with applications, Appl. Math. Comput. 219 (2012) 656-667.
\bibitem{hu}
L. Huang, W. So, On left eigenvalues of a quaternionic matrix, Linear Algebra Appl.
323 (2001) 105-116.

 \bibitem{so}
W. So,  Quaternionic left eigenvalue problem, Southeast
Asian Bulletin of Mathematics 29 (2005)  555-565.
 \bibitem{wo}
R. M. W. Wood, Quaternionic eigenvalues, Bull. Lond. Math. Soc. 17
(1985)137-138.

\bibitem{br}
J.L. Brenner, Matrices of quaternions, Pac. J. Math. 1 (1951) 329-335.
 \bibitem{ma} E. Mac\'{i}as-Virg\'{o}s, M.J. Pereira-S\'{a}ez,
A topological approach to left eigenvalues of quaternionic matrices, Linear Multilinear Algebra, (2013) DOI:10.1080/03081087.2012.753599.
\bibitem{ba} A. Baker, Right eigenvalues for quaternionic matrices: a topological ap-
proach, Linear Algebra Appl. 286 (1999) 303-309.
 \bibitem{dra} T. Dray, C. A. Manogue, The octonionic eigenvalue problem, Advances
in Applied Clifford Algebras 8(2) (1998) 341-364.
\bibitem{zh} F. Zhang, Quaternions and matrices of quaternions, Linear Algebra
Appl. 251 (1997) 21-57.
\bibitem{far} D.R. Farenick, B.A.F. Pidkowich, The spectral theorem in quaternions, Linear
Algebra Appl. 371 (2003)
75-102.
\bibitem{fa} F. O. Farid,   Q.W. Wang,  F. Zhang, On the eigenvalues of quaternion matrices,
Linear Multilinear Algebra, 59(4) (2011)  451- 473.
  \bibitem{la}P. Lancaster, M. Tismenitsky,  Theory of matrices, Acad. Press., New York, 1969.
\bibitem{ca}
Carl D. Meyer Jr.,  Limits and the index of a square
matrix, SIAM J. Appl. Math. 26(3) (1974) 506-515.








\end{thebibliography}
\end{document}